\documentclass[draft,preprint,12pt,numbers,sort&compress]{elsarticle}
\usepackage{mathrsfs}
\usepackage{amssymb}
\usepackage{amsthm}
\usepackage{mathrsfs}
\usepackage[centertags]{amsmath}
\usepackage{amsfonts}

\usepackage{color}

\input amssym.def
\input amssym.tex

\date{}

\newtheorem{Theorem}{Theorem}[section]

\newtheorem{Lemma}{Lemma}[section]

\newcommand\R{\mbox{\bf R}}
\newcommand\N{\mbox{\bf N}}
\newcommand\SR{\mbox{\scriptsize\bf R}}

\newcommand{\definition}{{\lower .5ex
  \hbox{$\>\>\stackrel{\triangle}{=}\>\>$} }}
\newcommand\supp{\mathop{\rm supp}}

\begin{document}

\baselineskip=22pt
\thispagestyle{empty}

\mbox{}
\bigskip

\begin{center}{\Large\bf Sharp  well-posedness  and ill-posedness of the }\\[1ex]
{\Large \bf Cauchy problem for the higher-order KdV}\\[1ex]
{\Large \bf  equations on the circle}

{Wei YAN$^a$,\quad Minjie JIANG$^a$, \quad
\quad Yongsheng LI$^b$ and Jianhua HUANG$^c$}\\[2ex]

{$^a$School of Mathematics and Information Science, Henan Normal University,}\\
{Xinxiang, Henan 453007, P. R. China}\\[1ex]

{$^b$Department of Mathematics, South China University of Technology,}\\
{Guangzhou, Guangdong 510640, P. R. China}\\[1ex]

{$^c$College of Science, National University of Defense and Technology,}\\
{ Changsha, Hunan 410073, P. R. China}\\[1ex]

\end{center}

\bigskip
\bigskip

\noindent{\bf Abstract.}  In this paper, we investigate the Cauchy problem
for the higher-order KdV-type equation
\begin{eqnarray*}
      u_{t}+(-1)^{j+1}\partial_{x}^{2j+1}u
     + \frac{1}{2}\partial_{x}(u^{2})
     = 0,j\in N^{+},x\in\mathbf{T}= [0,2\pi \lambda)
\end{eqnarray*}
 with low regularity data and $\lambda\geq 1$. Firstly, we show that the Cauchy
 problem for the periodic higher-order KdV equation
 is locally well-posed in $H^{s}(\mathbf{T})$ with $s\geq -j+\frac{1}{2},j\geq2.$
 By using some new Strichartz estimate and  some new function spaces,
 we also show that the Cauchy problem for the periodic higher-order KdV
  equation is ill-posed in $H^{s}(\mathbf{T})$ with $s<-j+\frac{1}{2},j\geq2$
  in the sense that  the solution map  is   $C^{3}.$ The result
  of this paper improves  the result of  \cite{H} with $j\geq2$.
\bigskip


\bigskip
\noindent {\bf Short Title:} Cauchy problem for higher-order KdV-type equation

\bigskip
\noindent {\bf Corresponding Author:} W. YAN

\bigskip
\noindent {\bf Email Address:}yanwei19821115@sina.cn

\bigskip
\noindent{\bf Fax number:} +86-0373-3326174

\bigskip

\leftskip 0 true cm \rightskip 0 true cm

\newpage{}

\begin{center}{\Large\bf Sharp  well-posedness  and ill-posedness of the }\\[1ex]
{\Large \bf Cauchy problem for the higher-order KdV}\\[1ex]
{\Large \bf  equations on the circle}

{Wei YAN$^a$,\quad Minjie JIANG$^a$, \quad
\quad Yongsheng LI$^b$ and Jianhua HUANG$^c$}\\[2ex]

{$^a$School of Mathematics and Information Science, Henan Normal University,}\\
{Xinxiang, Henan 453007, P. R. China}\\[1ex]

{$^b$Department of Mathematics, South China University of Technology,}\\
{Guangzhou, Guangdong 510640, P. R. China}\\[1ex]

{$^c$College of Science, National University of Defense and Technology,}\\
{ Changsha, Hunan 410073, P. R. China}\\[1ex]

\end{center}

\noindent{\bf Abstract.}
 In this paper, we investigate the Cauchy problem
for the higher-order KdV-type equation
\begin{eqnarray*}
      u_{t}+(-1)^{j+1}\partial_{x}^{2j+1}u
     + \frac{1}{2}\partial_{x}(u^{2})
     = 0,j\in N^{+},x\in [0,2\pi \lambda)
\end{eqnarray*}
 with low regularity data. Firstly, we show that the Cauchy
 problem for the periodic higher-order KdV equation
 is locally well-posed in $H^{s}(\mathbf{T})$ with
 $s\geq -j+\frac{1}{2},j\geq2.$
 By using some new Strichartz estimate and  some new function spaces,
 we also  show that the Cauchy problem for the
 periodic higher-order KdV
  equation is ill-posed in $H^{s}(\mathbf{T})$ with
  $s<-j+\frac{1}{2},j\geq2$ in the sense that  the solution map  is   $C^{3}.$ The result
  of this paper improves  the result of  \cite{H} with $j\geq2$.

\bigskip

{\large\bf 1. Introduction}
\bigskip

\setcounter{Theorem}{0} \setcounter{Lemma}{0}

\setcounter{section}{1}

In this paper, we consider the Cauchy problem for the periodic
higher-order KdV type equation
\begin{eqnarray}
&&u_{t}+(-1)^{j+1}\partial_{x}^{2j+1}u
     + \frac{1}{2}\partial_{x}(u^{2})
     = 0,\label{1.01}\\
    &&u(x,0)=u_{0}(x),\quad x\in \mathbf{T}=[0,2\pi\lambda), \label{1.02}
\end{eqnarray}
where $j\geq 2,j\in N$ and $\lambda\geq1.$
When $j=1,$
equation (\ref{1.01})  reduces to the Korteweg-de Vries (KdV) equation
\begin{eqnarray}
u_{t}+\partial_{x}^{3}u+\frac{1}{2}\partial_{x}(u^{2})=0.\label{1.03}
\end{eqnarray}
KdV equation possesses the bi-Hamiltonian structure and completely integrable,
thus, it possesses infinite conservation laws.

In recent some years, many people have paid more attention to
the Cauchy problem for the KdV
equation, for instance, see \cite{Bourgain93, Bourgain97, KPV1996,KPV2001,CKSTT,H,Kis,T} and the references therein.
 Using the Fourier
 restriction norm method introduced in \cite{B,Bourgain93} by Bourgain,
Kenig et. al. \cite{KPV1996}  proved that
   the Cauchy problem for the
KdV equation  is locally  well-posed in $H^{s}$ with $s>-\frac{3}{4}$ on the  real line  and  the
Cauchy problem for the periodic KdV equation is
locally well-posed in $H^{s}(\mathbf{T})$ with
$s\geq-\frac{1}{2}.$    Bourgain \cite{Bourgain97}
proved that the Cauchy problem for the periodic KdV
  equation is ill-posed in $H^{s}([0,2\pi))$
with $s< -\frac{1}{2}.$ By using the I-method,
  Colliander et.al. \cite{CKSTT} proved that
 the Cauchy problem for the
periodic KdV  equation is globally well-posed
 in $H^{s}(\mathbf{T})$ with $s\geq -\frac{1}{2}.$ Recently,
Kappeler and  Topalov \cite{KT2006} proved that
   the  Cauchy  problem for the KdV  equation  is locally
well-posed  in $H^{s}(\mathbf{T})$   with $s\geq -1$.
Guo \cite{G} and Kishimoto \cite{Kis} proved that
 the Cauchy problem for the KdV  equation is
 globally well-posed in $H^{-3/4}(\R).$
Molinet \cite{Molinet} proved that the  Cauchy
 problem for the KdV  equation  is locally  well-posed  in
$H^{s}(\mathbf{T})$   with $s\geq -1$ and
  ill-posed in $H^{s}(\mathbf{T})$   with $s< -1$.

Using the  Fourier restriction norm method, Hirayama \cite{H}
  proved that (\ref{1.01}) is locally
well-posed in $H^{s}(\mathbf{T})$ with $s\geq-\frac{j}{2}.$
In this paper, as in \cite{BT,Kato}, combining the  new function spaces introduced in this paper,
  the Strichartz estimate established in this paper with  the fixed point Theorem, we show that
the Cauchy problem for (\ref{1.01}) is
 locally well-posed in $H^{s}(\mathbf{T})$ with
$s\geq-j+\frac{1}{2}$ with $j\geq 2,j\in N$;
 we also show that the Cauchy problem for
 (\ref{1.01}) is ill-posed in $H^{s}(\mathbf{T})$
 with $s<-j+\frac{1}{2}$ with $j\geq 2,j\in N$
 in the sense that  the solution map  is  $C^{3}.$

We give some notations before presenting the main results. $C>0$
 may vary from line to line.  $0<\epsilon<\frac{1}{100j}.$ $A\sim B$ denotes that $|B|\leq |A|\leq 4|B|$.
 $A\gg B$ denotes that $|A|\geq 4|B|.$ $a\vee b={\rm max}\left\{a,b\right\}.$
 $a\wedge b={\rm min}\left\{a,b\right\}.$
 Throughout this paper,
 $\dot{Z}:=Z- \{ 0\}$ and $\dot{Z}^{+}:=Z^{+}- \{ 0\}$.
  Denote by
 $(dk)_{\lambda}$ the normalized counting measure on $\dot{Z_{\lambda}}=
 \frac{\dot{Z}}{\lambda}$:
 \begin{eqnarray*}
 \int a(k)(dk)_{\lambda}=\frac{1}{\lambda}\sum_{k\in  \dot{Z}_{\lambda}}a(k).
 \end{eqnarray*}
 Let
 \begin{eqnarray*}
 \mathscr{F}_{x}f(k)=\int_{0}^{2\pi\lambda}e^{- i kx}f(x)dx.
 \end{eqnarray*}
denote the Fourier transformation of a function $f$ on $[0,2\pi\lambda).$
 Let
 \begin{eqnarray*}
 f(x)=\int e^{ i kx} \mathscr{F}_{x}f(k)(dk)_{\lambda}=\frac{1}{\lambda}\sum_{k \in
 \dot{Z}_{\lambda}}e^{ i kx}\mathscr{F}_{x}f(k).
 \end{eqnarray*}
 Let
 \begin{eqnarray*}
 \mathscr{F}_{t}f(\tau)=\int_{\SR}e^{- i t\tau}f(t)dt
 \end{eqnarray*}
 denote
 the Fourier transformation of a function $f$
  with the respect to the time variable.
   Let
 \begin{eqnarray*}
 f(t)=\int e^{ i t\tau} \mathscr{F}_{t}f(\tau)d\tau.
 \end{eqnarray*}
Let
\begin{eqnarray*}
S(t)\phi(x)=\int e^{i kx }e^{i(-1)^{j+1}t k^{2j+1}}\mathscr{F}_{x}\phi(k)(dk)_{\lambda}.
\end{eqnarray*}
We define the space-time  Fourier transform
$\mathscr{F}f(k,\tau)$ for $k\in \dot{Z}$ and $\tau\in \R$ by
\begin{eqnarray*}
\mathscr{F}f(k,\tau)=\int\int_{0}^{2\pi\lambda}e^{-i kx}e^{-i \tau t}f(x,t)dxdt
\end{eqnarray*}
and this transformation is inverted by
\begin{eqnarray*}
f(x,t)=\int\int e^{i kx}e^{i \tau t}\mathscr{F}f(k,\tau)(dk)_{\lambda}d\tau.
\end{eqnarray*}
It is easily checked that
\begin{eqnarray*}
&&\|f\|_{L^{2}(\mathbf{T})}=\|\mathscr{F}_{x}f\|_{L^{2}((dk)_{\lambda})},\\
&&\int_{0}^{2\pi \lambda}f(x)\overline{g(x)}dx=\int \mathscr{F}_{x}f(k)\overline{\mathscr{F}_{x}
f(k)}(dk)_{\lambda},\\
&&\mathscr{F}_{x}(fg)=\mathscr{F}_{x}f*\mathscr{F}_{x}g=\int \mathscr{F}_{x}f(k-k_{1})
\mathscr{F}_{x}g(k_{1})(dk_{1})_{\lambda}.
\end{eqnarray*}
Let
\begin{eqnarray*}
&&P(k)=k^{2j+1},\sigma=\tau-P(k),\quad \sigma_{j}=\tau_{j}-P(k_{j}),\\
&&D_{1}=\left\{(\tau,k)\in \R\times \dot{Z}:|\tau-P(k)|\leq \frac{2(2j+1)}{3}|k|^{2j},|k|\geq1\right\},\\
&&D_{2}=\left\{(\tau,k)\in \R\times \dot{Z}: \frac{2(2j+1)}{3}|k|^{2j}< |\tau-P(k)|\leq 2(2j+1)|k|^{2j+1},|k|\geq1\right\},\\
&&D_{3}=\left\{(\tau,k)\in \R\times \dot{Z}:|\tau-P(k)|> 2(2j+1)|k|^{2j+1},|k|\geq1\right\},\\
&&D_{4}=\left\{(\tau,k)\in \R\times \dot{Z}:|\tau-P(k)|> 2(2j+1)|k|^{2j+1},\frac{1}{\lambda}\leq |k|\leq1\right\},\\
&&D_{5}=\left\{(\tau,k)\in \R\times \dot{Z}:|\tau-P(k)|\leq 2(2j+1)|k|^{2j+1},\frac{1}{\lambda}\leq |k|\leq1\right\},\\
&&\mathscr{F}(\Lambda^{-1})f=\langle \sigma \rangle ^{-1}\mathscr{F}f.
\end{eqnarray*}
We define the Sobolev space $H^{s}(\mathbf{T})$ with the norm
\begin{eqnarray*}
\|f\|_{H^{s}( \mathbf{T})}=\|\langle k\rangle^{s}\mathscr{F}_{x}f(k)\|_{L^{2}(dk)_{\lambda}}
\end{eqnarray*}
and define the $X_{s,b}$ spaces for $2\pi\lambda$-periodic KdV equipped with the norm
\begin{eqnarray*}
\|u\|_{X_{s,b}(\mathbf{T} \times \SR)}=\left\|\langle k\rangle^{s} \left\langle \tau
-P(k)\right\rangle^{b}\mathscr{F}u(k,\tau)\right\|_{L^{2}(dk)_{\lambda}(d\tau)}.
\end{eqnarray*}
We define the $Z^{s}$ space equipped with the norm
\begin{eqnarray*}
\|u\|_{Z^{s}}=\|P_{D_{1}\cup D_{5}}u\|_{X_{s,\frac{2j-1}{2j}}}+\|P_{D_{2}}u\|_{X_{(1-2j)s-1,s+1}}+
\|P_{D_{3}\cup D_{4}}u\|_{X_{-\frac{s}{j}-1,\frac{s}{j}+1}}+\left\|u\right\|_{Y^{s}},
\end{eqnarray*}
where $j\geq 2$ and $\|u\|_{Y^{s}}=\left\|\langle k\rangle^{s}
\mathscr{F}u(k,\tau)\right\|_{L^{2}(k)L^{1}(\tau)}.$

\noindent We define $Z^{s}([0,T])$ by the following norm:
\begin{eqnarray*}
\|u\|_{Z^{s}([0,T])}:={\rm inf}
\left\{\|v\|_{Z^{s}}:\qquad u(t)=v(t) \qquad on \qquad t\in I\right\}.
\end{eqnarray*}

The main result of this paper are as follow.

\begin{Theorem}\label{Thm1}
Let $s\geq -j+\frac{1}{2},j\geq2$  and $u_{0}$ be
$2\pi\lambda$-periodic
function.
Then the Cauchy problems (\ref{1.01})(\ref{1.02})
are locally well-posed in $H^{s}(\mathbf{T})$.
\end{Theorem}

\begin{Theorem}\label{Thm2}
Let $s< -j+\frac{1}{2},j\geq2$  and
$u_{0}$ be $2\pi $-periodic function.
Then the solution map $S_{t}$:  of the
Cauchy problem for (\ref{1.01}) is not $C^{3}$
at zero. More precisely,
for any $T > 0,$ the solution map
$S_{t} : u_{0}\in  H^{s}(0,2\pi )\rightarrow u
\in C([0,T]; H^{s}(0,2\pi ))$
is not $C^{3}$ at zero.
\end{Theorem}

\noindent {\bf Remark:} Theorem 1.1 is sharp
in the sense of Theorem 1.2.

The rest of the paper is arranged as follows.
In Section 2,  we present some
preliminaries. In Section 3, we establish three
important  bilinear estimates. In Section 4,
we present the proof of Theorem 1.1.
In Section 5, we present the proof of Theorem 1.2.

\bigskip
\bigskip
\bigskip

 \noindent{\large\bf 2. Preliminaries }

\setcounter{equation}{0}

\setcounter{Theorem}{0}

\setcounter{Lemma}{0}

\setcounter{section}{2}

In this section, we give some preliminaries which palys a paramount role
 in establishing Lemmas 3.1, 3.2.

\begin{Lemma}\label{Lemma2.1}
Let $u_{l}$ with $l=1,2$ be $L^{2}(\dot{Z}\times \R)$-real valued functions.
Then for any $(l_{1},l_{2}) \in \N^{2}$, we derive
\begin{eqnarray}
      \left\|(\Psi_{l_{1}}u_{1})*(\Psi_{l_{2}}u_{2})\right\|_{L_{k\tau}^{2}}
      \leq C\left(2^{l_{1}}
      \wedge2^{l_{2}}\right)^{1/2}\left(2^{l_{1}}\vee2^{l_{2}}\right)^{\frac{1}{2(2j+1)}}
      \|\Psi_{l_{1}}u_{1}\|_{L^{2}}\|\Psi_{l_{2}}u_{2}\|_{L^{2}}.
        \label{2.01}
\end{eqnarray}
\end{Lemma}

For the proof of  Lemma 2.1, we refer the readers to Lemma 2.1 of \cite{LYLH}.
\begin{Lemma}\label{Lemma2.2}
Let $u(x,t),v(x,t)$ be  $2\pi $-periodic functions and $a+b\geq \frac{j+1}{2j+1}$
and ${\rm min}\{a,b\}>\frac{1}{2(2j+1)}$. Then, we derive
\begin{eqnarray}
     && \left\|uv\right\|_{L_{xt}^{2}}\leq C\|u\|_{X_{0,a}([0,2\pi) \times \SR)}
     \|v\|_{X_{0,b}([0,2\pi) \times \SR)},
       \label{2.02}\\
       &&\left\|uv\right\|_{X_{0,-a}}\leq C\|u\|_{X_{0,b}([0,2\pi) \times \SR)}
       \|v\|_{L_{xt}^{2}}.
       \label{2.03}
\end{eqnarray}
\end{Lemma}
{\bf Proof.} From Lemma 2.1, we derive that
\begin{eqnarray}
&&\|uv\|_{L_{xt}^{2}}\leq C \sum_{l_{1}\geq 0}\sum_{l_{2}\geq 0}
\left\|(\Psi_{l_{1}}u_{1})*(\Psi_{l_{2}}u_{2})\right\|_{L_{xt}^{2}}\nonumber\\
&&\leq C\sum_{l_{1}\geq 0}\sum_{l_{2}\geq 0}\left(2^{l_{1}}
      \wedge2^{l_{2}}\right)^{1/2}\left(2^{l_{1}}\vee2^{l_{2}}\right)^{\frac{1}{2(2j+1)}}
      \|\Psi_{l_{1}}u_{1}\|_{L^{2}}\|\Psi_{l_{2}}u_{2}\|_{L^{2}}\label{2.04}.
\end{eqnarray}
Let $M_{j}=2^{l_{j}}$  with  $j=1,2.$ Without loss of generality, we can assume that $M_{1}\geq M_{2}$
and $M_{1}=NM_{2}$ and $u_{M_{j}}=\Psi_{l_{j}}u_{j}$ with $j=1,2.$
\begin{eqnarray}
&&\|uv\|_{L_{xt}^{2}}\leq  C\sum_{M_{1},M_{2}\geq 1}M_{1}^{\frac{1}{2(2j+1)}}
      M_{2}^{1/2}\|u_{M_{1}}\|_{L^{2}}\|v_{M_{2}}\|_{L^{2}}\nonumber\\
      &&\leq C\sum_{N,M_{2}\geq 1}M_{2}^{\frac{j+1}{2j+1}}N^{\frac{1}{2(2j+1)}}
      \|u_{NM_{2}}\|_{L^{2}}\|v_{M_{2}}\|_{L^{2}}\nonumber\\
      &&\leq C\sum_{M_{2},N\geq 1}N^{\frac{1}{2(2j+1)}-a}M_{2}^{\frac{j+1}{2j+1}-a-b}
      (NM_{2})^{a}\|u_{NM_{2}}\|_{L^{2}}M_{2}^{b}\|v_{M_{2}}\|_{L^{2}}\nonumber\\
      &&\leq C\|u\|_{X_{0,a}([0,2\pi)\times \SR)}\|u\|_{X_{0,b}( [0,2\pi)\times \SR)}.\label{2.05}
\end{eqnarray}
We can have (\ref{2.03})  by duality.

We have completed the proof of Lemma 2.2.

\begin{Lemma}\label{Lemma2.3}
Let $u(x,t),v(x,t)$ be  $2\pi\lambda $-periodic functions and $a+b\geq \frac{j+1}{2j+1}$
and ${\rm min}\{a,b\}>\frac{1}{2(2j+1)}$. Then, we obtain
\begin{eqnarray}
     && \left\|uv\right\|_{L_{xt}^{2}}\leq C\|u\|_{X_{0,a}(\mathbf{T} \times \SR)}
     \|v\|_{X_{0,b}(\mathbf{T}  \times \SR)},
       \label{2.02}\\
       &&\left\|uv\right\|_{X_{0,-a}}\leq C\|u\|_{X_{0,b}(\mathbf{T}  \times \SR)}
       \|v\|_{L_{xt}^{2}}.
       \label{2.03}
\end{eqnarray}
\end{Lemma}

By using a similar technique of Lemma 3.4 in \cite{Molinet} and
Lemma 2.2, we  have Lemma 2.3.

\begin{Lemma}\label{Lemma2.4}Let $s\in \R$ and  $T>0$ and $u(x,t),v(x,t)$ be  $2\pi\lambda $-periodic functions. Then, we have
\begin{eqnarray*}
\left\| \int_{0}^{t}S(t-\tau)\partial_{x}(uv)d\tau\right\|_{Z^{s}([0,T])}
\leq C\|\partial_{x}\Lambda^{-1}(uv)\|_{Z^{s}([0,T])}.
\end{eqnarray*}
\end{Lemma}

For the proof of Lemma 2.4, we refer the readers to \cite{BT}.

\begin{Lemma}\label{Lemma2.5}Let $s\in \R$ and $j\geq 2,j\in Z$ and $u(x,t)$ be  $2\pi\lambda $-periodic function. Then, we derive
\begin{eqnarray}
&&\|u\|_{X^{s,\frac{1}{2j}}}\leq C\|u\|_{Z^{s}}\leq C\|u\|_{X^{s,\frac{2j-1}{2j}}},\label{2.06}\\
&&\|u\|_{X^{s,\frac{1}{2}}(D_{1}\bigcup D_{2})}\leq C\|u\|_{Z^{s}(D_{1}\bigcup D_{2})}.\label{2.07}
\end{eqnarray}
\end{Lemma}
{\bf Proof.} We firstly have that (\ref{2.06}). When $\supp \mathscr{F}u\subset D_{1}$,
since $\frac{2j-1}{2j}\geq \frac{1}{2j},$ we derive that
$\|u\|_{X_{s,\frac{2j-1}{2j}}}\geq \|u\|_{X_{s,\frac{1}{2j}}}$. When $\supp \mathscr{F}u\subset D_{2}$,
 since $-j+\frac{1}{2}\leq s\leq-\frac{j}{2}$,  we derive that
$\langle \sigma\rangle ^{s+\frac{2j-1}{2j}}\geq C\langle k\rangle^{2js+1}$ which yields
that $\langle k\rangle^{s}\langle \sigma\rangle^{\frac{1}{2j}}\leq C\langle k\rangle^{(1-2j)s-1}\langle \sigma\rangle^{s+1}$,
thus, we derive that $\|u\|_{X_{(1-2j)s-1,s+1}}\geq \|u\|_{X_{s,\frac{1}{2j}}}$. When $\supp \mathscr{F}u\subset D_{3}$,
since $-j+\frac{1}{2}\leq s\leq-\frac{j}{2}$,  we obtain that
$\langle \sigma\rangle ^{\frac{s}{j}+\frac{2j-1}{2j}}\geq C\langle k\rangle^{s+1+\frac{s}{j}}$ which yields
that $|k|^{s}\langle \sigma\rangle^{\frac{1}{2j}}\leq C\langle k\rangle^{-\frac{s}{j}-1}\langle \sigma\rangle^{\frac{s}{j}+1}$,
thus, we have that $\|u\|_{X_{-\frac{s}{j}-1,\frac{s}{j}+1}}\geq \|u\|_{X_{s,\frac{1}{2j}}}$.
Consequently, we have that $\|u\|_{Z^{s}}\geq C\|u\|_{X_{s,\frac{1}{2j}}}.$
When $\supp \mathscr{F}u\subset D_{2}$,
 since $-j+\frac{1}{2}\leq s\leq-\frac{j}{2}$,  we obtain that
$\langle \sigma\rangle ^{s+\frac{2j-1}{2j}}\geq C\langle k\rangle^{2js+1}$ which yields
that $\langle k\rangle^{(1-2j)s-1}\langle \sigma\rangle^{s+1}\leq C\langle k\rangle^{s}
\langle \sigma\rangle^{\frac{2j-1}{2j}}$,
thus, we have that $\|u\|_{X_{(1-2j)s-1,s+1}}\leq C\|u\|_{X_{s,\frac{2j-1}{2j}}}$.
When $\supp \mathscr{F}u\subset D_{3}$,
since $-j+\frac{1}{2}\leq s\leq-\frac{j}{2}$,  we derive that
$\langle \sigma\rangle ^{-\frac{s}{j}-\frac{1}{2j}}\geq C\langle k\rangle^{-s-1-\frac{s}{j}}$ which yields
that $\langle k\rangle^{-\frac{s}{j}-1}\langle \sigma\rangle^{\frac{s}{j}+1}
\leq C\langle k\rangle^{s}\langle \sigma\rangle^{\frac{2j-1}{2j}}$,
thus, we obtain that $\|u\|_{X_{-\frac{s}{j}-1,\frac{s}{j}+1}}\leq C \|u\|_{X_{s,\frac{2j-1}{2j}}}$.
Consequently, we derive that $\|u\|_{Z^{s}}\leq C\|u\|_{X_{s,\frac{2j-1}{2j}}}.$
By Cauchy-Schwartz inequality with respect to $\tau$, we derive that
$\|\langle k\rangle^{s}\mathscr{F}u\|_{l_{k}^{2}l_{\tau}^{1}}\leq C
\|u\|_{X_{s,\frac{2j-1}{2j}}},$  consequently, we have that $\|u\|_{Z^{s}}
\leq C\|u\|_{X^{s,\frac{2j-1}{2j}}}.$  Now we show (\ref{2.07}). When
$\supp \mathscr{F}u\subset D_{1}$, since $\frac{2j-1}{2j}\geq \frac{1}{2},$ we derive that
$\|u\|_{X_{s,\frac{2j-1}{2j}}}\geq \|u\|_{X_{s,\frac{1}{2}}}$. When
$\supp \mathscr{F}u\subset D_{2}$, since $s\geq -j+\frac{1}{2}$, we derive that
$\langle k\rangle^{s}\langle\sigma\rangle^{1/2}\leq C\langle k\rangle
^{(1-2j)s-1}\langle\sigma\rangle^{s+1},$ consequently, we derive that
$\|u\|_{X_{(1-2j)s,\>s+1}}\geq \|u\|_{X_{s,\frac{1}{2}}}$.

We have completed the proof of Lemma 2.5.

{\bf Remark:}
The conclusion of case $j=2$ of (\ref{2.06})-(\ref{2.07})  can be found in \cite{Kato}.
\begin{Lemma}\label{Lemma2.6}
Assume that  $s\in \R$ and  $T>0$. Then, we have that
\begin{eqnarray*}
\left\| S(t)\phi\right\|_{Z^{s}([0,T])}\leq C\|\phi\|_{H^{s}(\mathbf{T})}.
\end{eqnarray*}
\end{Lemma}
{\bf Proof.} From the definition of  Lemma 2.5, we have that $X^{\frac{2j-1}{2j}}([0,T])
\hookrightarrow Z^{s}([0,T])\hookrightarrow C([0,T]:H^{s}(\mathbf{T})).$

We have completed the proof of Lemma 2.6.

\begin{Lemma}\label{Lemma2.7}
Let
\begin{eqnarray*}
\sigma=\tau-k^{2j+1},\sigma_{1}= \tau_{1}-k_{1}^{2j+1},\sigma_{2}= \tau_{1}-k_{2}^{2j+1}.
\end{eqnarray*}
 Then, we have
\begin{eqnarray*}
3{\rm max}\left\{|\sigma|,|\sigma_{1}|,|\sigma_{2}|\right\}\geq|\sigma-\sigma_{1}-\sigma_{2}|=
\left|k^{2j+1}-k_{1}^{2j+1}-k_{2}^{2j+1}\right|\geq (2j+1)|kk_{1}^{j}k_{2}^{j}|.
\end{eqnarray*}
\end{Lemma}

For the proof  of Lemma 2.7, we refer the readers to Lemma 2.4 of \cite{LYY}.
From Lemma 2.7, we have that one of three following cases must occur:
\begin{eqnarray}
&& (a): |\sigma|={\rm max}\left\{|\sigma|,|\sigma_{1}|,|\sigma_{2}|\right\},\label{2.010}\\
&& (b): |\sigma_{1}|={\rm max}\left\{|\sigma|,|\sigma_{1}|,|\sigma_{2}|\right\},\label{2.011}\\
&& (c): |\sigma_{2}|={\rm max}\left\{|\sigma|,|\sigma_{1}|,|\sigma_{2}|\right\}.\label{2.012}
\end{eqnarray}

\bigskip
\bigskip

\noindent{\large\bf 3. Bilinear estimates }

\setcounter{equation}{0}

 \setcounter{Theorem}{0}

\setcounter{Lemma}{0}

 \setcounter{section}{3}
 This section is devoted to establishing some significant bilinear estimates which are used to derive the Theorem 1.1.
 \begin{Lemma}\label{Lemma3.1}
Let $j\geq 2$ and $-j+\frac{1}{2}\leq s\leq-\frac{j}{2}$ and $u_{j}(x,t)(j=1,2)$ be  $2\pi\lambda $-periodic functions. Then, we have
\begin{eqnarray}
      \left\|\Lambda^{-1}\partial_{x}(\prod_{j=1}^{2}u_{j})\right\|_{X^{s}}\leq C\prod\limits_{j=1}^{2}\|u_{j}\|_{Z^{s}},
        \label{3.01}
\end{eqnarray}
here $C>0$,  which  is  independent of  $\lambda$, $\left\|\cdot\right\|_{X^{s}}$ is the norm removing $\left\|\cdot\right\|_{Y^{s}}$ from $\left\|\cdot\right\|_{Z^{s}}.$
\end{Lemma}
{\bf Proof.} Obviously, $\left(\R\times\dot{Z}_{\lambda}\right)^{2}\subset \bigcup\limits_{j=1}^{8}\Omega_{j},$
where
\begin{eqnarray*}
&&\Omega_{1}=\left\{(\tau_{1},k_{1},\tau,k)\in \left(\R\times\dot{Z_{\lambda}}\right)^{2}:{\rm max}\left\{|k_{1}|, |k|\right\}\leq1\right\},\\
&&\Omega_{2}=\left\{(\tau_{1},k_{1},\tau,k)\in \left(\R\times\dot{Z_{\lambda}}\right)^{2}\cap \Omega_{1}^{c}:|k_{1}|\sim |k_{2}|\gg |k|\geq1\right\},\\
&&\Omega_{3}=\left\{(\tau_{1},k_{1},\tau,k)\in \left(\R\times\dot{Z_{\lambda}}\right)^{2}\cap \Omega_{1}^{c}:|k_{1}|\sim |k_{2}|\gg |k|,1\geq|k|\geq\frac{1}{\lambda}\right\},\\
&&\Omega_{4}=\left\{(\tau_{1},k_{1},\tau,k)\in \left(\R\times\dot{Z_{\lambda}}\right)^{2}\cap \Omega_{1}^{c}:|k|\sim |k_{2}|\gg |k_{1}|\geq 1\right\},\\
&&\Omega_{5}=\left\{(\tau_{1},k_{1},\tau,k)\in \left(\R\times\dot{Z_{\lambda}}\right)^{2}\cap \Omega_{1}^{c}:|k|\sim |k_{2}|\gg |k_{1}|,1\geq|k_{1}|\geq \frac{1}{\lambda}\right\},\\
&&\Omega_{6}=\left\{(\tau_{1},k_{1},\tau,k)\in \left(\R\times\dot{Z_{\lambda}}\right)^{2}\cap \Omega_{1}^{c}:|k|\sim |k_{1}|\gg |k_{2}|\geq 1\right\},\\
&&\Omega_{7}=\left\{(\tau_{1},k_{1},\tau,k)\in \left(\R\times\dot{Z_{\lambda}}\right)^{2}\cap \Omega_{1}^{c}:|k|\sim |k_{1}|\gg |k_{2}|,1\geq|k_{2}|\geq \frac{1}{\lambda}\right\},\\
&&\Omega_{8}=\left\{(\tau_{1},k_{1},\tau,k)\in \left(\R\times\dot{Z_{\lambda}}\right)^{2}\cap
\Omega_{1}^{c}:|k|\sim |k_{1}|\sim |k_{2}|\geq 1\right\}.
\end{eqnarray*}
(1) In region $\Omega_{1}$.
By using Lemma 2.5,  from the definition of $Z^{s},$ since
${\rm max}\left\{|k_{1}|,|k|\right\}\leq 1$ and
$-j+\frac{1}{2}\leq s\leq-\frac{j}{2}$, we have that
\begin{eqnarray*}
&&\left\|\Lambda^{-1}\partial_{x}(\prod_{j=1}^{2}u_{j})\right\|_{X^{s}}\leq C\left\|(1-\partial_{x}^{2})^{-\frac{1}{2}}\partial_{x}
(\prod_{j=1}^{2}u_{j})\right\|_{X_{s,\frac{2j-1}{2j}}}\nonumber\\
&&\leq C\left\|k\langle \sigma\rangle^{-\frac{1}{2j}}
\left(\mathscr{F}u_{1}*\mathscr{F}u_{2}\right)
\right\|_{l_{k}^{2}L_{\tau}^{2}}\nonumber\\
&&\leq C\||k|\|_{l_{k}^{2}}\left\|\mathscr{F}u_{1}*\mathscr{F}u_{2}\right\|_{l_{k}^{\infty}L_{\tau}^{2}}\nonumber\\
&&\leq C\|\mathscr{F}u_{1}\|_{l_{k}^{2}L_{\tau}^{2}}\|\mathscr{F}u_{2}\|_{l_{k}^{2}L_{\tau}^{1}}\nonumber\\
&&\leq C\|u_{1}\|_{X_{s,\frac{1}{2j}}}\|u_{2}\|_{Y^{s}}\leq C\prod_{j=1}^{2}\|u_{j}\|_{Z^{s}}.
\end{eqnarray*}
(2) In region $\Omega_{2}$. In this  region, we consider  (a)-(c) of Lemma 2.7, respectively.

\noindent (a) Case $|\sigma|={\rm max}\left\{|\sigma|,|\sigma_{1}|,|\sigma_{2}|\right\},$ we have
$\supp \left[\mathscr{F}u_{1}*\mathscr{F}u_{2}\right]\subset D_{3}.$

\noindent
When $\supp \mathscr{F}u_{j} \subset \Omega_{1}\cup \Omega_{2}$ with $j=1,2$,
 by using Lemmas 2.5, 2.7, 2.3,  since $-j+\frac{1}{2}\leq s\leq-\frac{j}{2}$, we have that
\begin{eqnarray*}
&&\left\|\Lambda^{-1}\partial_{x}(\prod_{j=1}^{2}u_{j})\right\|_{X^{s}}\leq C\left\|\langle k\rangle^{-\frac{s}{j}}\langle\sigma\rangle^{\frac{s}{j}}
\left[\mathscr{F}u_{1}*\mathscr{F}u_{2}\right]\right\|_{l_{k}^{2}L_{\tau}^{2}}
\nonumber\\
&&\leq C\|(J^{s}u_{1})(J^{s}u_{2})\|_{L_{xt}^{2}}\nonumber\\
&&\leq C\|u_{1}\|_{X_{s,\frac{1}{2}}}\|u_{2}\|_{X_{s,\frac{1}{2(2j+1)}+\epsilon}}\leq C\|u_{1}\|_{X_{s,\frac{1}{2}}}\|u_{2}\|_{X_{s,\frac{1}{2j}}}\leq C\prod_{j=1}^{2}\|u_{j}\|_{Z^{s}}.
\end{eqnarray*}
When $\supp \mathscr{F}u_{1} \subset \Omega_{3}$,
 by using Lemmas 2.5, 2.7 and the  Young inequality, since $-j+\frac{1}{2}\leq s\leq-\frac{j}{2}$, we have that
\begin{eqnarray*}
&&\left\|\Lambda^{-1}\partial_{x}(\prod_{j=1}^{2}u_{j})\right\|_{X^{s}}\leq C\left\|\langle k\rangle^{-\frac{s}{j}}\langle\sigma\rangle^{\frac{s}{j}}
(\mathscr{F}u_{1}*\mathscr{F}u_{2})\right\|_{l_{k}^{2}L_{\tau}^{2}}\nonumber\\
&&\leq C\left\|\mathscr{F}u_{1}*\left[\langle k\rangle^{2s}\mathscr{F}u_{2}\right]\right\|_{l_{k}^{2}L_{\tau}^{2}}\nonumber\\
&&\leq C\left\|\left[\langle k\rangle^{-\frac{s}{j}-1}\langle\sigma\rangle^{\frac{s}{j}+1}u_{1}\right]*\left[\langle k\rangle^{-2j}\mathscr{F}u_{2}\right]\right\|_{l_{k}^{2}L_{\tau}^{2}}
\nonumber\\
&&\leq C\|u\|_{X_{-\frac{s}{j}-1,\frac{s}{j}+1}}\|\langle k\rangle^{-2j}\mathscr{F}u_{2}\|_{l_{k}^{1}L_{\tau}^{1}}\nonumber\\
&&\leq C\|u\|_{X_{-\frac{s}{j}-1,\frac{s}{j}+1}}\|\langle k\rangle^{-2j-s}\langle k\rangle^{s}\mathscr{F}u_{2}\|_{l_{k}^{1}L_{\tau}^{1}}\nonumber\\
&&\leq C\|u\|_{X_{-\frac{s}{j}-1,\frac{s}{j}+1}}\|\langle k\rangle^{s}\mathscr{F}u_{2}\|_{l_{k}^{2}L_{\tau}^{1}}\nonumber\\
&&\leq
C\prod_{j=1}^{2}\|u_{j}\|_{Z^{s}}.
\end{eqnarray*}
When $\supp \mathscr{F}u_{2} \subset \Omega_{3}$,
this case can be proved similarly to $\supp \mathscr{F}u_{1} \subset \Omega_{3}$.

\noindent (b) Case $|\sigma_{1}|={\rm max}\left\{|\sigma|,|\sigma_{1}|,|\sigma_{2}|\right\},$
in this case, we consider the following cases
\begin{eqnarray*}
(i): |\sigma_{1}|>4{\rm max}\left\{|\sigma|,|\sigma_{2}|\right\},
(ii):|\sigma_{1}|\leq4{\rm max}\left\{|\sigma|,|\sigma_{2}|\right\},
\end{eqnarray*}
respectively.

\noindent
When (i) occurs:
if $\supp u_{1}\subset D_{1}$ which yields that $1\leq |k|\leq C$, by using Lemmas 2.5,
2.7, 2.3, since $ -j+\frac{1}{2} \leq s\leq-\frac{j}{2} ,$
 we have that
\begin{eqnarray*}
&&\left\|\Lambda^{-1}\partial_{x}(\prod_{j=1}^{2}u_{j})\right\|_{X^{s}}\leq C\left\|\langle k\rangle ^{s+1}\langle \sigma \rangle ^{-\frac{1}{2j}}
(\mathscr{F}u_{1}*\mathscr{F}u_{2})\right\|_{l_{k}^{2}L_{\tau}^{2}}\nonumber\\&&
\leq C\left\|\langle \sigma \rangle ^{-\frac{1}{2j}}(\langle k\rangle^{s}\langle\sigma\rangle^{\frac{2j-1}{2j}}\mathscr{F}u_{1})*(\langle k\rangle ^{-s-2j+1}\mathscr{F}u_{2})\right\|_{l_{k}^{2}L_{\tau}^{2}}\nonumber\\
&&\leq C\left\|\left(J^{s}\Lambda ^{\frac{2j-1}{2j}}u_{1}\right)
\left(J^{-s-2j+1}u_{2}\right)\right\|_{X_{0,-\frac{1}{2j}}}\nonumber\\
&&\leq C\|u_{1}\|_{X_{s,\frac{2j-1}{2j}}}\|u_{2}\|_{X_{-s-2j+1,\frac{1}{2}}}\nonumber\\
&&\leq C\|u_{1}\|_{X_{s,\frac{2j-1}{2j}}}\|u_{2}\|_{X_{s,\frac{1}{2}}}\leq
C\prod_{j=1}^{2}\|u_{j}\|_{Z^{s}};
\end{eqnarray*}
if $\supp u_{1}\subset D_{2},$  by using Lemmas 2.5, 2.7, 2.3, since $-j+\frac{1}{2}\leq s\leq-\frac{j}{2}\leq -1$,
 we have that
\begin{eqnarray*}
&&\left\|\Lambda^{-1}\partial_{x}(\prod_{j=1}^{2}u_{j})\right\|_{X^{s}}\leq C\left\|\langle k\rangle ^{s+1}\langle \sigma \rangle ^{-\frac{1}{2j}}
\left[\mathscr{F}u_{1}*\mathscr{F}u_{2}\right]\right\|_{l_{k}^{2}L_{\tau}^{2}}\nonumber\\&&
\leq C\left\|\left(J^{(1-2j)s-1}\Lambda ^{s+1}u_{1}\right)\left(J^{-s-2j+1}u_{2}\right)
\right\|_{X_{0,-\frac{1}{2j}}}\nonumber\\
&&\leq C\|u_{1}\|_{X_{(1-2j)s-1,s+1}}\|u_{2}\|_{X_{-s-2j+1,\frac{1}{2}}}\nonumber\\
&&\leq C\|u_{1}\|_{X_{(1-2j)s-1,s+1}}\|u_{2}\|_{X_{s,\frac{1}{2}}}\leq
C\prod_{j=1}^{2}\|u_{j}\|_{Z^{s}}.
\end{eqnarray*}
When (ii) occurs: we have $|\sigma_{1}|\sim |\sigma|$ or $|\sigma_{1}|\sim |\sigma_{2}|$.

\noindent
When $|\sigma_{1}|\sim |\sigma|$  is valid, this case can be proved similarly to
$|\sigma|={\rm max}\left\{|\sigma|,|\sigma_{1}|,|\sigma_{2}|\right\}.$
When $|\sigma_{1}|\sim |\sigma_{2}|$, if
$\supp u_{1}\subset D_{1}$ which yields that $1\leq |k|\leq C$,
by using Lemmas 2.5, 2.7, 2.3, since $-j+\frac{1}{2}\leq s\leq -\frac{j}{2},$
 we have that
\begin{eqnarray*}
&&\left\|\Lambda^{-1}\partial_{x}(\prod_{j=1}^{2}u_{j})\right\|_{X^{s}}\leq C
\left\|\Lambda^{-1}\partial_{x}(\prod_{j=1}^{2}u_{j})\right\|_{X^{s}}\leq C
\left\|\langle k\rangle ^{s+1}\langle \sigma \rangle ^{-\frac{1}{2j}}
(\mathscr{F}u_{1}*\mathscr{F}u_{2})\right\|_{l_{k}^{2}L_{\tau}^{2}}\nonumber\\&&
\leq C\left\|(\langle k\rangle^{s}\langle\sigma\rangle^{\frac{2j-1}{2j}}\mathscr{F}u_{1})*(\langle k\rangle ^{-s-2j+1}\mathscr{F}u_{2})\right\|_{X_{0,-\frac{1}{2j}}}\nonumber\\
&&\leq C\left\|\left(J^{s}\Lambda ^{\frac{2j-1}{2j}}u_{1}\right)
\left(J^{-s-2j+1}u_{2}\right)\right\|_{X_{0,-\frac{1}{2j}}}\nonumber\\
&&\leq C\|u_{1}\|_{X_{(1-2j)s-1,s+1}}\|u_{2}\|_{X_{-s-2j+1,\frac{1}{2}}}\nonumber\\
&&\leq C\|u_{1}\|_{X_{s,\frac{2j-1}{2j}}}\|u_{2}\|_{X_{s,\frac{1}{2}}}\leq
C\prod_{j=1}^{2}\|u_{j}\|_{Z^{s}};
\end{eqnarray*}
if $\supp u_{1}\subset D_{2},$  by using Lemmas 2.5, 2.7, 2.3, since $-j+\frac{1}{2}\leq s\leq-\frac{j}{2}\leq -1$,
 we have that
\begin{eqnarray*}
&&\left\|\Lambda^{-1}\partial_{x}(\prod_{j=1}^{2}u_{j})\right\|_{X^{s}}\leq C
\left\|\langle k\rangle ^{s+1}\langle \sigma \rangle ^{-\frac{1}{2j}}
\left[\mathscr{F}u_{1}*\mathscr{F}u_{2}\right]\right\|_{l_{k}^{2}L_{\tau}^{2}}\nonumber\\&&
\leq C\left\|\left(J^{(1-2j)s-1}\Lambda ^{s+1}u_{1}\right)\left(J^{-s-2j+1}u_{2}\right)
\right\|_{X_{0,-\frac{1}{2j}}}\nonumber\\
&&\leq C\|u_{1}\|_{X_{(1-2j)s-1,s+1}}\|u_{2}\|_{X_{-s-2j+1,\frac{1}{2}}}\nonumber\\
&&\leq C\|u_{1}\|_{X_{(1-2j)s-1,s+1}}\|u_{2}\|_{X_{s,\frac{1}{2}}}\nonumber\\
&&\leq
C\prod_{j=1}^{2}\|u_{j}\|_{Z^{s}};
\end{eqnarray*}
if $\mathscr{F}u_{1}\subset D_{3}$, then  $\mathscr{F}u_{2}\subset D_{3}$,
by using the H\"older  inequality  and the  Young inequality,
 from Lemma 2.5, since $-j+\frac{1}{2}\leq s\leq-\frac{j}{2},$ we  have that
\begin{eqnarray*}
&&\left\|\Lambda^{-1}\partial_{x}(\prod_{j=1}^{2}u_{j})\right\|_{X^{s}}\leq C\left\|\langle k\rangle ^{s+1}\langle \sigma\rangle ^{-\frac{1}{2j}}
\left[\mathscr{F}u_{1}*\mathscr{F}u_{2}\right]\right\|_{l_{k}^{2}L_{\tau}^{2}}\nonumber\\
&&\leq C\left\|\langle k\rangle ^{s+\frac{3}{2}+\epsilon}\langle \sigma \rangle^{-\frac{1}{2j}+\frac{1}{2}+\epsilon}
\left[\mathscr{F}u_{1}*\mathscr{F}u_{2}\right]
\right\|_{l_{k}^{\infty}L_{\tau}^{\infty}}\nonumber\\
&&\leq C\left\|\left(\langle k\rangle ^{s+j+1-\frac{1}{2j}+(2j+2)\epsilon}\mathscr{F}u_{1}
\right)*\mathscr{F}u_{2}\right\|_{l_{k}^{\infty}l_{\tau}^{\infty}}\nonumber\\&&\leq C
\left\|\langle k\rangle ^{-3s+1-\frac{1}{2j}-3j+(2j+2)\epsilon}\right\|_{l_{k}^{\infty}}
\prod_{j=1}^{2}\|u_{j}\|_{X_{(1-2j)s-1,s+1}}\nonumber\\&&\leq C
\prod_{j=1}^{2}\|u_{j}\|_{X_{(1-2j)s-1,s+1}}\leq C\prod_{j=1}^{2}\|u_{j}\|_{Z^{s}}.
\end{eqnarray*}
(c) Case $|\sigma_{2}|={\rm max}\left\{|\sigma|,|\sigma_{1}|,|\sigma_{2}|\right\}.$  This case can be proved similarly to case (b).

\noindent (3) Region $\Omega_{3}$.
We  consider $|k|\leq |k_{1}|^{-2j}$  and  $|k_{1}|^{-2j}<|k|\leq 1,$ respectively.

\noindent When $|k|\leq |k_{1}|^{-2j}$, by using Lemmas 2.5,  2.7, since $-j+\frac{1}{2}\leq s\leq-\frac{j}{2},$  we have that
\begin{eqnarray*}
&&\left\|\Lambda^{-1}\partial_{x}(\prod_{j=1}^{2}u_{j})\right\|_{X^{s}}\leq C\left\||k|\langle \sigma\rangle ^{-\frac{1}{2j}}
\left[\mathscr{F}u_{1}*\mathscr{F}u_{2}\right]\right\|_{l_{k}^{2}L_{\tau}^{2}}\nonumber\\
&&\leq C\left\|\left[\langle k\rangle^{-(2j-1)}\mathscr{F}u_{1}*\langle k\rangle^{-(2j-1)}\mathscr{F}u_{2}\right]
\right\|_{l_{k}^{\infty}L_{\tau}^{2}}\nonumber\\
&&\leq C\|u_{1}\|_{X_{1-2j,0}}\|u_{2}\|_{Y^{1-2j}}
\leq C\|u_{1}\|_{X_{1-2j,0}}\|u_{2}\|_{Y^{s}};
\end{eqnarray*}
if $\supp \mathscr{F}u_{1}\subset D_{1},$ since $ -j+\frac{1}{2}\leq s\leq-\frac{j}{2},$ then $\|u\|_{X_{1-2j,0}}\leq C\|u\|_{X_{s,\frac{1}{2}}};$
if $\supp \mathscr{F}u_{1}\subset D_{2},$ since $-j+\frac{1}{2}\leq s\leq-\frac{j}{2};$
 then $\|u\|_{X_{1-2j,0}}\leq C\|u\|_{X_{(1-2j)s-1,s+1}};$
if $\supp \mathscr{F}u_{1}\subset D_{3},$ since $-j+\frac{1}{2}\leq s\leq-\frac{j}{2};$
then $\|u\|_{X_{1-2j,0}}\leq C\|u\|_{X_{-\frac{s}{j}-1,\frac{s}{j}+1}},$
thus, according to the definition of $Z^{s},$  we have that
\begin{eqnarray*}
\|u_{1}\|_{X_{1-2j,0}}\|u_{2}\|_{Y^{s}}\leq C\prod_{j=1}^{2}\|u_{j}\|_{Z^{s}}.
\end{eqnarray*}
Now we consider the case $|k_{1}|^{-2j}<|k|\leq 1.$ In this case,
we consider cases (a)-(c) of Lemma 2.7, respectively.

\noindent When   (a) occurs:   in this case
$\supp \left[\mathscr{F}u_{1}*\mathscr{F}u_{2}\right]\subset D_{4}$,
by using the H\"older  inequality and the
Young inequality and Lemma 2.7,  since
$|k|\leq 1$ and $\frac{s}{j}+1\leq 0$
and $-j+\frac{1}{2}\leq s\leq-\frac{j}{2},$  we have that
\begin{eqnarray*}
&&\left\|\Lambda^{-1}\partial_{x}
(\prod_{j=1}^{2}u_{j})\right\|_{X^{s}}
\leq C\left\||k|\langle k\rangle ^{-\frac{s}{j}-1}\langle \sigma \rangle ^{\frac{s}{j}}\left[\mathscr{F}u_{1}*\mathscr{F}u_{2}\right]\right\|_{l_{k}^{2}L_{\tau}^{2}}\nonumber\\
&&\leq C\left\||k|^{1+\frac{s}{j}} \left[(|k|^{s}\mathscr{F}u_{1})*(|k|^{s}\mathscr{F}u_{2})\right]
\right\|_{l_{k}^{2}L_{\tau}^{2}}\nonumber\\
&&\leq C\left\|\left[(|k|^{s}\mathscr{F}u_{1})*(|k|^{s}\mathscr{F}u_{2})\right]
\right\|_{l_{k}^{\infty}L_{\tau}^{2}}\leq C\|u_{1}\|_{X_{s,0}}\|u_{2}\|_{Y^{s}};
\end{eqnarray*}
if $\supp \mathscr{F}u_{1}\subset D_{1},$ since $-j+\frac{1}{2}\leq s\leq-\frac{j}{2},$ then $\|u\|_{X_{s,0}}\leq C\|u\|_{X_{s,\frac{1}{2}}};$
if $\supp \mathscr{F}u_{1}\subset D_{2},$ since $-j+\frac{1}{2}\leq s\leq-\frac{j}{2},$
 then $\|u\|_{X_{s,0}}\leq C\|u\|_{X_{(1-2j)s-1,s+1}};$
if $\supp \mathscr{F}u_{1}\subset D_{3},$ since $-j+\frac{1}{2}\leq s\leq-\frac{j}{2},$
then $\|u\|_{X_{s,0}}\leq C\|u\|_{X_{-\frac{s}{j}-1,\frac{s}{j}+1}};$
thus, according to the definition of $Z^{s},$  we have that
\begin{eqnarray*}
\|u_{1}\|_{X_{s,0}}\|u_{2}\|_{Y^{s}}\leq C\prod_{j=1}^{2}\|u_{j}\|_{Z^{s}}.
\end{eqnarray*}
When (b)  occurs:  we   consider the case
$|\sigma_{1}|>4{\rm max}\left\{|\sigma|,|\sigma_{2}|\right\}$  and
$|\sigma_{1}|\leq4{\rm max}\left\{|\sigma|,|\sigma_{2}|\right\}$.
When $|\sigma_{1}|>4{\rm max}\left\{|\sigma|,|\sigma_{2}|\right\}$ is valid, we have that $\supp \mathscr{F}u_{j} \subset D_{1}$  with $j=1,2.$

\noindent
If $\supp \left[\mathscr{F}u_{1}*\mathscr{F}u_{2}\right]\subset D_{4}$, by
by using the H\"older  inequality and the Young inequality,  since $|k|\leq 1$ and $1+\frac{s}{j}\geq 0,$  by using Lemma 2.7, since $-j+\frac{1}{2}\leq s\leq-\frac{j}{2},$ we have that
\begin{eqnarray*}
&&\left\|\Lambda^{-1}\partial_{x}(\prod_{j=1}^{2}u_{j})\right\|_{X^{s}}\leq C\left\||k|\langle k\rangle ^{-\frac{s}{j}-1}\langle \sigma \rangle ^{\frac{s}{j}}
\left[\mathscr{F}u_{1}*\mathscr{F}u_{2}\right]
\right\|_{l_{k}^{2}L_{\tau}^{2}}\nonumber\\
&&\leq C\left\||k|^{1+\frac{s}{j}} \left[(|k|^{s}\mathscr{F}u_{1})*(|k|^{s}\mathscr{F}u_{2})\right]
\right\|_{l_{k}^{2}L_{\tau}^{2}}\nonumber\\
&&\leq C\left\|\left[(|k|^{s}\mathscr{F}u_{1})*(|k|^{s}\mathscr{F}u_{2})\right]
\right\|_{l_{k}^{\infty}L_{\tau}^{2}}\nonumber\\
&&\leq C\|u_{1}\|_{X_{s,0}}\|u_{2}\|_{Y^{s}}\nonumber\\
&&\leq C\prod_{j=1}^{2}\|u_{j}\|_{X_{s,\frac{2j-1}{2j}}}\leq C\prod_{j=1}^{2}\|u_{j}\|_{Z^{s}};
\end{eqnarray*}
if $\supp \left[\mathscr{F}u_{1}*\mathscr{F}u_{2}\right]\subset D_{5}$, by using Lemma
 2.3, since $-j+\frac{1}{2}\leq s\leq-\frac{j}{2},$ we have that
\begin{eqnarray*}
&&\left\|\Lambda^{-1}\partial_{x}(\prod_{j=1}^{2}u_{j})\right\|_{X^{s}}
\leq C\left\||k|\langle k\rangle ^{s}\langle \sigma \rangle^{-\frac{1}{2j}}
\left[\mathscr{F}u_{1}*\mathscr{F}u_{2}\right]
\right\|_{l_{k}^{2}L_{\tau}^{2}}\nonumber\\
&&\leq C\left\|\langle \sigma \rangle^{-\frac{1}{2j}}
\left[(\langle k\rangle ^{s}\langle \sigma \rangle
^{\frac{2j-1}{2j}}\mathscr{F}u_{1})*
\left(\langle k\rangle^{-s-2j+1}\mathscr{F}u_{2}\right)
\right]\right\|_{l_{k}^{2}L_{\tau}^{2}}\nonumber\\
&&\leq C\left\|\left(J^{s}\Lambda^{\frac{2j-1}{2j}}u_{1}
\right)\left(J^{-s-2j+1}u_{2}\right)\right\|_{X_{0,-\frac{1}{2j}}}\nonumber\\
&&\leq C\|J^{s}\Lambda^{\frac{2j-1}{2j}}u_{1}\|_{L_{xt}^{2}}
\|J^{-s-2j+1}u_{2}\|_{X_{0,\frac{1}{2}}}\nonumber\\
&&\leq C\|u_{1}\|_{X_{s,\frac{2j-1}{2j}}}\|u_{2}\|_{X_{s,\frac{1}{2}}}\nonumber\\
&&\leq C\prod_{j=1}^{2}\|u_{j}\|_{X_{s,\frac{2j-1}{2j}}}\leq C\prod_{j=1}^{2}\|u_{j}\|_{Z^{s}}.
\end{eqnarray*}
When $|\sigma_{1}|\leq4{\rm max}\left\{|\sigma|,|\sigma_{2}|\right\}$ is valid, we have
 $|\sigma_{1}|\sim |\sigma|$ or  $|\sigma_{1}|\sim |\sigma_{2}|$.

\noindent If $|\sigma_{1}|\sim|\sigma|$, then this case an be proved similarly to
 case $|\sigma|={\rm max}\left\{|\sigma|,|\sigma_{1}|,|\sigma_{2}|\right\}.$

\noindent
  When $|\sigma_{1}|\sim|\sigma_{2}|$, we consider $\supp\mathscr{F}u_{1}\subset \Omega_{1}$
and $\supp\mathscr{F}u_{1}\subset \Omega_{2}$, respectively.

\noindent When $\supp\emph{}\mathscr{F}u_{1}\subset \Omega_{1}$,
if $\supp \left[\mathscr{F}u_{1}*\mathscr{F}u_{2}\right]\subset D_{5},$ since $-j+\frac{1}{2}\leq s\leq-\frac{j}{2},$
 we have that
\begin{eqnarray*}
&&\left\|\Lambda^{-1}\partial_{x}(\prod_{j=1}^{2}u_{j})\right\|_{X^{s}}
\leq C\left\||k|\langle k\rangle^{s}\langle \sigma \rangle^{-\frac{1}{2j}}\left[\mathscr{F}u_{1}*\mathscr{F}u_{2}\right]
\right\|_{l_{k}^{2}L_{\tau}^{2}}\nonumber\\
&&\leq C\left\|\left(J^{s}\Lambda^{\frac{2j-1}{2j}}u_{1}\right)\left(J^{-s-2j+1}u_{2}\right)
\right\|_{X_{0,-\frac{1}{2j}}}\nonumber\\
&&\leq C\|u_{1}\|_{X_{s,\frac{2j-1}{2j}}}\|u_{2}\|_{X_{-s-2j+1,\frac{1}{2}}}\nonumber\\
&&\leq C\|u_{1}\|_{X_{s,\frac{2j-1}{2j}}}\|u_{2}\|_{X_{s,\frac{1}{2}}}\nonumber\\&&
\leq C\prod_{j=1}^{2}\|u_{j}\|_{Z^{s}};
\end{eqnarray*}
if $\supp \left[\mathscr{F}u_{1}*\mathscr{F}u_{2}\right]\subset D_{4},$ since $-j+\frac{1}{2}\leq s\leq-\frac{j}{2}$,
 we have that
\begin{eqnarray*}
&&\left\|\Lambda^{-1}\partial_{x}(\prod_{j=1}^{2}u_{j})\right\|_{X^{s}}\leq C\left\||k|\langle k\rangle^{\frac{s}{j}}\langle \sigma \rangle^{\frac{s}{j}}\left[\mathscr{F}u_{1}*\mathscr{F}u_{2}\right]\right\|_{l_{k}^{2}L_{\tau}^{2}}\nonumber\\
&&\leq C\left\||k|\langle \sigma \rangle^{-\frac{1}{2}+\epsilon}\left[\mathscr{F}u_{1}*\mathscr{F}u_{2}\right]\right\|_{l_{k}^{2}L_{\tau}^{2}}\nonumber\\
&&\leq C\left\|\left(J^{s}\Lambda^{\frac{2j-1}{2j}}u_{1}\right)
\left(J^{-s-2j+1}u_{2}\right)\right\|_{X_{0,-\frac{1}{2}+\epsilon}}\nonumber\\
&&\leq C\|u_{1}\|_{X_{s,\frac{2j-1}{2j}}}\|u_{2}\|_{X_{-s-2j+1,\frac{2}{2j+1}}}\leq C\prod_{j=1}^{2}\|u_{j}\|_{Z^{s}}.
\end{eqnarray*}
\noindent When $\supp\mathscr{F}u_{1}\subset \Omega_{2}$,
if $\supp \left[\mathscr{F}u_{1}*\mathscr{F}u_{2}\right]\subset D_{5},$  since $-j+\frac{1}{2}\leq s\leq-\frac{j}{2},$
 we have that
\begin{eqnarray*}
&&\left\|\Lambda^{-1}\partial_{x}(\prod_{j=1}^{2}u_{j})\right\|_{X^{s}}
\leq C\left\||k|\langle k\rangle^{s}\langle \sigma \rangle^{-\frac{1}{2j}}\left[\mathscr{F}u_{1}*\mathscr{F}u_{2}\right]
\right\|_{l_{k}^{2}L_{\tau}^{2}}\nonumber\\
&&\leq C\left\||k|\langle k\rangle^{s}\langle \sigma \rangle^{-\frac{1}{2}+\epsilon}\left[\mathscr{F}u_{1}*\mathscr{F}u_{2}\right]
\right\|_{l_{k}^{2}L_{\tau}^{2}}\nonumber\\
&&\leq C\left\|\left(J^{(1-2j)s-1}\Lambda^{s+1}u_{1}\right)\left(J^{-s-2j+1}u_{2}\right)
\right\|_{X_{0,-\frac{1}{2}+\epsilon}}\nonumber\\
&&\leq C\|u_{1}\|_{X_{(1-2j)s-1,s+1}}\|u_{2}\|_{X_{-s-2j+1,\frac{2}{2j+1}}}\leq C\prod_{j=1}^{2}\|u_{j}\|_{Z^{s}};
\end{eqnarray*}
if $\supp \left[\mathscr{F}u_{1}*\mathscr{F}u_{2}\right]\subset D_{4},$
since $-j+\frac{1}{2}\leq s\leq-\frac{j}{2}$,
 we have that
\begin{eqnarray*}
&&\left\|\Lambda^{-1}\partial_{x}(\prod_{j=1}^{2}u_{j})\right\|_{X^{s}}
\leq C\left\||k|\langle k\rangle^{\frac{s}{j}}\langle \sigma \rangle^{\frac{s}{j}}\left[\mathscr{F}u_{1}*\mathscr{F}u_{2}\right]
\right\|_{l_{k}^{2}L_{\tau}^{2}}\nonumber\\
&&\leq C\left\||k|\langle \sigma \rangle^{-\frac{1}{2}+\epsilon}\left[\mathscr{F}u_{1}*\mathscr{F}u_{2}\right]
\right\|_{l_{k}^{2}L_{\tau}^{2}}\nonumber\\
&&\leq C\left\|\left(J^{(1-2j)s-1}\Lambda^{s+1}u_{1}\right)\left(J^{-s-2j+1}u_{2}\right)
\right\|_{X_{0,-\frac{1}{2}+\epsilon}}\nonumber\\
&&\leq C\|u_{1}\|_{X_{s,\frac{2j-1}{2j}}}\|u_{2}\|_{X_{-s-2j+1,\frac{2}{2j+1}}}\leq C\prod_{j=1}^{2}\|u_{j}\|_{Z^{s}}.
\end{eqnarray*}
(4) Region $\Omega_{4}$.
We consider   cases (a)-(c) of Lemma 2.7, respectively.

\noindent When (a) occurs: if $|\sigma|>4{\rm max}\left\{|\sigma_{1}|,|\sigma_{2}|\right\}$,
 then $\supp \mathscr{F}u_{j}\subset D_{1}\cup D_{2}$ with $j=1,2$  and  $\supp \left[\mathscr{F}u_{1}*\mathscr{F}u_{2}\right]\subset D_{2}.$
In this case, by using Lemma 2.3, since $-j+\frac{1}{2}\leq s\leq -\frac{j}{2},$ we have that
\begin{eqnarray*}
&&\left\|\Lambda^{-1}\partial_{x}(\prod_{j=1}^{2}u_{j})\right\|_{X^{s}}\leq C\left\|\langle k\rangle^{-(2j-1)s}\langle \sigma \rangle^{s}\left[\mathscr{F}u_{1}*\mathscr{F}u_{2}\right]\right\|_{l_{k}^{2}L_{\tau}^{2}}\nonumber\\
&&\leq C\left\|(J^{s}u_{1})(J^{s}u_{2})\right\|_{L_{xt}^{2}}\leq C\|u_{1}\|_{X_{s,\frac{1}{2j}}}\|u_{2}\|_{X_{s,\frac{1}{2}}}\leq C\prod_{j=1}^{2}\|u_{j}\|_{Z^{s}}.
\end{eqnarray*}
When $|\sigma|\leq4{\rm max}\left\{|\sigma_{1}|,|\sigma_{2}|\right\}$, we have  that $|\sigma|\sim |\sigma_{1}|$ or $|\sigma|\sim |\sigma_{2}|$.

\noindent
 When  $|\sigma|\sim |\sigma_{1}|$, if $\supp \left[\mathscr{F}u_{1}*\mathscr{F}u_{2}\right]\subset D_{2},$
 by using Lemma 2.3,  since $-j+\frac{1}{2}\leq s\leq -\frac{j}{2},$ we have that
\begin{eqnarray*}
&&\left\|\Lambda^{-1}\partial_{x}(\prod_{j=1}^{2}u_{j})\right\|_{X^{s}}\leq C\left\|\langle k\rangle^{-(2j-1)s}\langle \sigma \rangle^{s}
\left[\mathscr{F}u_{1}*\mathscr{F}u_{2}\right]\right\|_{l_{k}^{2}L_{\tau}^{2}} \nonumber\\&&\leq\left\|(J^{s}u_{1})(J^{s}u_{2})\right\|_{L_{xt}^{2}}\nonumber\\&&\leq C\|u_{1}\|_{X_{s,\frac{1}{2j}}}\|u_{2}\|_{X_{s,\frac{1}{2}}}\leq C\prod_{j=1}^{2}\|u_{j}\|_{Z^{s}};
\end{eqnarray*}
if $\supp \left[\mathscr{F}u_{1}*\mathscr{F}u_{2}\right]\subset D_{3},$
 by using Lemma 2.3,  since $-j+\frac{1}{2}\leq s\leq -\frac{j}{2},$ we have that
\begin{eqnarray*}
&&\left\|\Lambda^{-1}\partial_{x}(\prod_{j=1}^{2}u_{j})\right\|_{X^{s}}\leq C\left\|\langle k\rangle^{-\frac{s}{j}}\langle \sigma \rangle^{\frac{s}{j}}
\left[\mathscr{F}u_{1}*\mathscr{F}u_{2}\right]\right\|_{l_{k}^{2}L_{\tau}^{2}} \nonumber\\&&\leq\left\|(J^{-\frac{s}{j}-1}
\Lambda^{\frac{s}{j}+1}u_{1})(J^{-\frac{s}{j}-(2j-1)}u_{2})\right\|_{X_{0,-\frac{1}{2j}}}\nonumber\\&&\leq C\|u_{1}\|_{X_{-\frac{s}{j}-1,\frac{s}{j}+1}}\|u_{2}\|_{X_{s,\frac{1}{2}}}\leq C\prod_{j=1}^{2}\|u_{j}\|_{Z^{s}}.
\end{eqnarray*}
When $|\sigma|\sim |\sigma_{2}|$,
 if $\supp \left[\mathscr{F}u_{1}*\mathscr{F}u_{2}\right]\subset D_{2},$
 by using Lemma 2.3,  since $-j+\frac{1}{2}\leq s\leq -\frac{j}{2},$ we have that
\begin{eqnarray*}
&&\left\|\Lambda^{-1}\partial_{x}(\prod_{j=1}^{2}u_{j})\right\|_{X^{s}}\leq C\left\|\langle k\rangle^{-(2j-1)s}\langle \sigma \rangle^{s}
\left[\mathscr{F}u_{1}*\mathscr{F}u_{2}\right]\right\|_{l_{k}^{2}L_{\tau}^{2}} \nonumber\\&&\leq\left\|(J^{s}u_{1})(J^{s}u_{2})\right\|_{L_{xt}^{2}}\nonumber\\&&\leq C\|u_{1}\|_{X_{s,\frac{1}{2j}}}\|u_{2}\|_{X_{s,\frac{1}{2}}}\leq C\prod_{j=1}^{2}\|u_{j}\|_{Z^{s}};
\end{eqnarray*}
if $\supp \left[\mathscr{F}u_{1}*\mathscr{F}u_{2}\right]\subset D_{3},$
 by using Lemma 2.3,  since $-j+\frac{1}{2}\leq s\leq -\frac{j}{2},$ we have that
\begin{eqnarray*}
&&\left\|\Lambda^{-1}\partial_{x}(\prod_{j=1}^{2}u_{j})\right\|_{X^{s}}\leq C\left\|\langle k\rangle^{-\frac{s}{j}}\langle \sigma \rangle^{\frac{s}{j}}
\left[\mathscr{F}u_{1}*\mathscr{F}u_{2}\right]\right\|_{l_{k}^{2}L_{\tau}^{2}} \nonumber\\&&\leq\left\|(J^{-\frac{s}{j}-(2j-1)}u_{1})(J^{-\frac{s}{j}-1}
\Lambda^{\frac{s}{j}+1}u_{2})\right\|_{X_{0,-\frac{1}{2j}}}\nonumber\\&&\leq C\|u_{1}\|_{X_{s,\frac{1}{2}}}\|u_{2}\|_{X_{-\frac{s}{j}-1,\frac{s}{j}+1}}\leq C\prod_{j=1}^{2}\|u_{j}\|_{Z^{s}}.
\end{eqnarray*}

\noindent (b): $|\sigma_{1}|={\rm max}\left\{|\sigma|,|\sigma_{1}|,|\sigma_{2}|\right\}.$
If $|\sigma_{1}|>4{\rm max}\left\{|\sigma|,|\sigma_{2}|\right\}$,
 then $\supp \mathscr{F}u_{1}\subset D_{3}$.
  We consider
  \begin{eqnarray*}
  \supp \left[\mathscr{F}u_{1}*\mathscr{F}u_{2}\right]\subset D_{1},\supp \left[\mathscr{F}u_{1}*\mathscr{F}u_{2}\right]\subset D_{2},
  \end{eqnarray*}
  respectively.

  \noindent
  When $\supp \left[\mathscr{F}u_{1}*\mathscr{F}u_{2}\right]\subset D_{1},$
 by using Lemma 2.3,   since $-j+\frac{1}{2}\leq s\leq-\frac{j}{2},$  we have that
\begin{eqnarray*}
&&\left\|\Lambda^{-1}\partial_{x}(\prod_{j=1}^{2}u_{j})\right\|_{X^{s}}\leq C\left\|\langle k\rangle^{s+1}\langle \sigma \rangle^{-\frac{1}{2j}}
\left[\mathscr{F}u_{1}*\mathscr{F}u_{2}\right]\right\|_{l_{k}^{2}L_{\tau}^{2}}
\nonumber\\
&&\leq C\left\|(J^{-\frac{s}{j}-1}\Lambda ^{\frac{s}{j}+1}u_{1})(J^{-s-2j+1}u_{2})\right\|_{X^{0,-\frac{1}{2j}}}
\nonumber\\&&\leq C\|u_{1}\|_{X_{\frac{s}{j}+1,\frac{s}{j}+1}}\|u_{2}\|_{X_{s,\frac{1}{2}}}
\leq C\prod_{j=1}^{2}\|u_{j}\|_{Z^{s}};
\end{eqnarray*}
if $\supp \left[\mathscr{F}u_{1}*\mathscr{F}u_{2}\right]\subset D_{2},$
 by using Lemma 2.3,   since $-j+\frac{1}{2}\leq s\leq-\frac{j}{2},$  we have that
\begin{eqnarray*}
&&\left\|\Lambda^{-1}\partial_{x}(\prod_{j=1}^{2}u_{j})\right\|_{X^{s}}\leq C\left\|\langle k\rangle^{s+1}\langle \sigma \rangle^{-\frac{1}{2j}}
\left[\mathscr{F}u_{1}*\mathscr{F}u_{2}\right]\right\|_{l_{k}^{2}L_{\tau}^{2}}
\nonumber\\
&&\leq C\left\|(J^{-\frac{s}{j}-1}\Lambda ^{\frac{s}{j}+1}u_{1})(J^{-s-2j+1}u_{2})\right\|_{X^{0,-\frac{1}{2j}}}
\nonumber\\&&\leq C\|u_{1}\|_{X_{\frac{s}{j}+1,\frac{s}{j}+1}}\|u_{2}\|_{X_{s,\frac{1}{2}}}
\leq C\prod_{j=1}^{2}\|u_{j}\|_{Z^{s}}.
\end{eqnarray*}
When $|\sigma_{1}|\leq4{\rm max}\left\{|\sigma|,|\sigma_{2}|\right\}$,
we have that $|\sigma_{1}|\sim |\sigma|$ or $|\sigma_{1}|\sim |\sigma_{2}|$.

\noindent When  $|\sigma_{1}|\sim |\sigma|$, this case can be proved similarly
 to case $|\sigma|={\rm max}\left\{|\sigma|,|\sigma_{1}|,|\sigma_{2}|\right\}.$

\noindent When $|\sigma_{1}|\sim |\sigma_{2}|$,
 if $\supp \left[\mathscr{F}u_{1}*\mathscr{F}u_{2}\right]\subset D_{1},$
 by using Lemma 2.3, since $-j+\frac{1}{2}\leq s\leq-\frac{j}{2},$ we have that
\begin{eqnarray*}
&&\left\|\Lambda^{-1}\partial_{x}(\prod_{j=1}^{2}u_{j})\right\|_{X^{s}}\leq C\left\|\langle k\rangle^{s+1}\langle \sigma \rangle^{-\frac{1}{2j}}
\left[\mathscr{F}u_{1}*\mathscr{F}u_{2}\right]\right\|_{l_{k}^{2}L_{\tau}^{2}}
 \nonumber\\&&\leq\left\|(J^{-\frac{s}{j}-1}
\Lambda^{\frac{s}{j}+1}u_{1})(J^{-s-2j+1}u_{2})\right\|_{X^{0,-\frac{1}{2j}}}
\leq C\|u_{1}\|_{X_{s,\frac{1}{2j}}}\|u_{2}\|_{X_{s,\frac{1}{2}}}\leq C\prod_{j=1}^{2}\|u_{j}\|_{Z^{s}};
\end{eqnarray*}
if $\supp \left[\mathscr{F}u_{1}*\mathscr{F}u_{2}\right]\subset D_{2},$
 by using Lemma 2.3, since $-j+\frac{1}{2}\leq s\leq-\frac{j}{2},$ we have that
\begin{eqnarray*}
&&\left\|\Lambda^{-1}\partial_{x}(\prod_{j=1}^{2}u_{j})\right\|_{X^{s}}\leq C\left\|\langle k\rangle^{-(2j-1)s}\langle \sigma \rangle^{s}
\left[\mathscr{F}u_{1}*\mathscr{F}u_{2}\right]\right\|_{l_{k}^{2}L_{\tau}^{2}} \nonumber\\&&\leq\left\|(J^{s}u_{1})(J^{s}u_{2})\right\|_{L_{xt}^{2}}\nonumber\\&&\leq C\|u_{1}\|_{X_{s,\frac{1}{2j}}}\|u_{2}\|_{X_{s,\frac{1}{2}}}\leq C\prod_{j=1}^{2}\|u_{j}\|_{Z^{s}};
\end{eqnarray*}
if $\supp \left[\mathscr{F}u_{1}*\mathscr{F}u_{2}\right]\subset D_{3},$
 by using Lemma 2.3, since $-j+\frac{1}{2}\leq s\leq-\frac{j}{2}$ and $\epsilon<\frac{1}{100j},$  we have that
\begin{eqnarray*}
&&\left\|\Lambda^{-1}\partial_{x}(\prod_{j=1}^{2}u_{j})\right\|_{X^{s}}\leq C\left\|\langle k\rangle ^{s+1}\langle \sigma\rangle ^{-\frac{1}{2j}}
\left[\mathscr{F}u_{1}*\mathscr{F}u_{2}\right]\right\|_{l_{k}^{2}L_{\tau}^{2}}\nonumber\\
&&\leq C\left\|\langle k\rangle ^{s+\frac{3}{2}+\epsilon}\langle \sigma \rangle^{-\frac{1}{2j}+\frac{1}{2}+\epsilon}
\left[\mathscr{F}u_{1}*\mathscr{F}u_{2}\right]
\right\|_{l_{k}^{\infty}L_{\tau}^{\infty}}\nonumber\\
&&\leq C\left\|\left(\langle k\rangle ^{s+j+1-\frac{1}{2j}+(2j+2)\epsilon}\mathscr{F}u_{1}
\right)*\mathscr{F}u_{2}\right\|_{l_{k}^{\infty}l_{\tau}^{\infty}}\nonumber\\&&\leq C
\left\|\langle k\rangle ^{-3s+1-\frac{1}{2j}-3j+(2j+2)\epsilon}\right\|_{l_{k}^{\infty}}
\prod_{j=1}^{2}\|u_{j}\|_{X_{-\frac{s}{j}-1,\frac{s}{j}+1}}\nonumber\\&&\leq C
\prod_{j=1}^{2}\|u_{j}\|_{X_{-\frac{s}{j}-1,\frac{s}{j}+1}}\leq C\prod_{j=1}^{2}\|u_{j}\|_{Z^{s}}.
\end{eqnarray*}
\noindent (c): $|\sigma_{2}|={\rm max}\left\{|\sigma|,|\sigma_{1}|,|\sigma_{2}|\right\}.$
In this case, we consider $|\sigma_{2}|\geq 4{\rm max}\left\{|\sigma|,|\sigma_{1}|\right\}$
and $|\sigma_{2}|< 4{\rm max}\left\{|\sigma|,|\sigma_{1}|\right\}$, respectively.

\noindent When $|\sigma_{2}|\geq 4{\rm max}\left\{|\sigma|,|\sigma_{1}|\right\}$,
$\mathscr{F}u_{2}\subset D_{2}$ and $\left[\mathscr{F}u_{1}*\mathscr{F}u_{2}\right]\subset D_{1}\bigcup D_{2}$, since $-j+\frac{1}{2}\leq s\leq-\frac{j}{2},$ we have
\begin{eqnarray*}
&&\left\|\Lambda^{-1}\partial_{x}(\prod_{j=1}^{2}u_{j})\right\|_{X^{s}}\leq C\left\|\langle k\rangle^{s+1}\langle \sigma\rangle ^{-\frac{1}{2j}}
\left[\mathscr{F}u_{1}*\mathscr{F}u_{2}\right]\right\|_{l_{k}^{2}L_{\tau}^{2}}\nonumber\\
&&\leq C\left\|\left[(J^{-s-2j+1}\mathscr{F}u_{1})*(\langle k\rangle^{(1-2j)s-1}\Lambda^{s+1}\mathscr{F}u_{2})\right]\right\|
_{X_{0,-\frac{1}{2j}}}\nonumber\\
&&\leq \|u_{1}\|_{X_{s,\frac{1}{2}}}\|u_{2}\|_{X^{(1-2j)s-1,s+1}}\nonumber\\
&&\leq C\prod_{j=1}^{2}\|u_{j}\|_{Z^{s}}.
\end{eqnarray*}
When $|\sigma_{2}|< 4{\rm max}\left\{|\sigma|,|\sigma_{1}|\right\}$, we have $|\sigma_{2}|\sim |\sigma_{1}|$ or $|\sigma_{2}|\sim |\sigma|.$

\noindent Case $|\sigma_{2}|\sim |\sigma_{1}|$  can be proved simialrly to $|\sigma_{1}|={\rm max}\left\{|\sigma|,|\sigma_{1}|,|\sigma_{2}|\right\}$.

\noindent
$|\sigma_{2}|\sim |\sigma|$ can be proved simialrly to $|\sigma|={\rm max}\left\{|\sigma|,|\sigma_{1}|,|\sigma_{2}|\right\}$.

\noindent(5) In region $\Omega_{5}$. In this case,
we consider the case $|k_{1}|\leq |k|^{-2j}$ and
$|k|^{-2j}< |k_{1}|\leq 1,$ respectively.

\noindent When $|k_{1}|\leq |k|^{-2j}$, by using the
Young inequality and Cauchy-Schwartz inequality as well as Lemmas 2.5 2.7, since $-j+\frac{1}{2}\leq s\leq-\frac{j}{2},$ we have that
\begin{eqnarray*}
&&\left\|\Lambda^{-1}\partial_{x}(\prod_{j=1}^{2}u_{j})\right\|_{X^{s}}\leq C\left\|\langle k\rangle^{s+1}\langle \sigma\rangle ^{-\frac{1}{2j}}
\left[\mathscr{F}u_{1}*\mathscr{F}u_{2}\right]\right\|_{l_{k}^{2}L_{\tau}^{2}}\nonumber\\
&&\leq C\left\|\left[(|k|^{-\frac{1}{2j}}\mathscr{F}u_{1})*(\langle k\rangle^{s}\mathscr{F}u_{2})\right]\right\|_{l_{k}^{2}L_{\tau}^{2}}\nonumber\\
&&\leq C\left\||k|^{-\frac{1}{2j}}\mathscr{F}u_{1}\right\|_{l_{k}^{1}l_{\tau}^{2}}\|u_{2}\|_{Y^{s}}\nonumber\\&&
\leq C\|\mathscr{F}u_{1}\|_{l_{n}^{2}L_{\tau}^{2}}\|u_{2}\|_{Y^{s}}
\leq C\prod_{j=1}^{2}\|u_{j}\|_{Z^{s}}.
\end{eqnarray*}
When  $|k|^{-2j}\leq |k_{1}|\leq 1,$
we consider   cases (a)-(c) of Lemma 2.7, respectively.

\noindent
When (a) occurs:
 by using the Young inequality and Cauchy-Schwartz inequality as well as Lemma 2.7, since $-j+\frac{1}{2}\leq s\leq-\frac{j}{2},$ we have that
\begin{eqnarray*}
&&\left\|\Lambda^{-1}\partial_{x}(\prod_{j=1}^{2}u_{j})\right\|_{X^{s}}\leq C\left\|\langle k\rangle^{s+1}\langle \sigma\rangle ^{-\frac{1}{2j}}
\left[\mathscr{F}u_{1}*\mathscr{F}u_{2}\right]\right\|_{l_{k}^{2}L_{\tau}^{2}}\nonumber\\
&&\leq C\left\|\left[(|k|^{-\frac{1}{2j}}\mathscr{F}u_{1})*(\langle k\rangle ^{s}\mathscr{F}u_{2})\right]
\right\|_{l_{k}^{2}L_{\tau}^{2}}\nonumber\\
&&\leq C\left\||k|^{-\frac{1}{2j}}\mathscr{F}u_{1}\right\|_{l_{k}^{1}l_{\tau}^{2}}\|u_{2}\|_{Y^{s}}\nonumber\\&&
\leq C\|\mathscr{F}u_{1}\|_{l_{k}^{2}L_{\tau}^{2}}\|u_{2}\|_{Y^{s}}
\leq C\prod_{j=1}^{2}\|u_{j}\|_{Z^{s}}.
\end{eqnarray*}
When (b) occurs:
 by using the Young inequality and Cauchy-Schwartz inequality and Lemma 2.7, since $-j+\frac{1}{2}\leq s\leq-\frac{j}{2},$ we have that
\begin{eqnarray*}
&&\left\|\Lambda^{-1}\partial_{x}(\prod_{j=1}^{2}u_{j})\right\|_{X^{s}}\leq C\left\|\langle k\rangle ^{s+1}\langle \sigma\rangle ^{-\frac{1}{2j}}
\left[\mathscr{F}u_{1}*\mathscr{F}u_{2}\right]\right\|_{l_{k}^{2}L_{\tau}^{2}}\nonumber\\
&&\leq C\left\|\left[(|k|^{-\frac{1}{2j}}\langle \sigma \rangle ^{\frac{1}{2j}}\mathscr{F}u_{1})*(\langle k\rangle ^{s}\mathscr{F}u_{2})\right]
\right\|_{l_{k}^{2}L_{\tau}^{2}}\nonumber\\
&&\leq C\left\||k|^{-\frac{1}{2j}}\langle \sigma \rangle ^{\frac{1}{2j}}\mathscr{F}u_{1}\right\|_{l_{k}^{1}l_{\tau}^{2}}\|u_{2}\|_{Y^{s}}\nonumber\\&&
\leq C\|\langle \sigma \rangle ^{\frac{1}{2j}}\mathscr{F}u_{1}\|_{l_{k}^{2}L_{\tau}^{2}}\|u_{2}\|_{Y^{s}}\nonumber\\&&
\leq C\|u_{1}\|_{X_{0,\frac{1}{2j}}}\|u_{2}\|_{Y^{s}}
\leq C\prod_{j=1}^{2}\|u_{j}\|_{Z^{s}}.
\end{eqnarray*}
When (c)  occurs: in this case $\langle k_{2}\rangle^{-s}\langle\sigma_{2}\rangle^{-\frac{1}{2j}}\leq |k_{1}|^{-\frac{1}{2j}}\langle k_{2}\rangle^{-s-1},$
by using the Young inequality and Cauchy-Schwartz inequality, since $-j+\frac{1}{2}\leq s\leq-\frac{j}{2},$ we have that
\begin{eqnarray*}
&&\left\|\Lambda^{-1}\partial_{x}(\prod_{j=1}^{2}u_{j})\right\|_{X^{s}}\leq C\left\|\langle k\rangle ^{s+1}\langle \sigma\rangle ^{-\frac{1}{2j}}
\left[\mathscr{F}u_{1}*\mathscr{F}u_{2}\right]\right\|_{l_{k}^{2}L_{\tau}^{2}}\nonumber\\
&&\leq C\left\|\left[(|k|^{-\frac{1}{2j}}\mathscr{F}u_{1})*(\langle k\rangle ^{s}\langle \sigma \rangle ^{\frac{1}{2j}}\mathscr{F}u_{2})\right]
\right\|_{l_{k}^{2}L_{\tau}^{2}}\nonumber\\
&&\leq C\left\||k|^{-\frac{1}{2j}}\mathscr{F}u_{1}\right\|_{l_{k}^{1}l_{\tau}^{1}}\|u_{2}\|_{X_{s,\frac{1}{2j}}}\nonumber\\&&
\leq C\|\mathscr{F}u_{1}\|_{l_{k}^{2}L_{\tau}^{1}}\|u_{2}\|_{X_{s,\frac{1}{2j}}}\nonumber\\&&
\leq C\|u_{1}\|_{Y^{s}}\|u_{2}\|_{X_{s,\frac{1}{2j}}}
\leq C\prod_{j=1}^{2}\|u_{j}\|_{Z^{s}}.
\end{eqnarray*}
(6)In region $\Omega_{6}$.
This region can be proved similarly to $\Omega_{4}$.

\noindent (7)In region $\Omega_{7}$.
This region can be proved similarly to $\Omega_{5}$.

\noindent (8)In region $\Omega_{8}$. We consider   cases (a)-(c) of Lemma 2.7, respectively.

\noindent When (a) occurs: $\supp\left( \mathscr{F}u_{1}*\mathscr{F}u_{2}\right) \subset D_{3}$.
If $|\sigma|>4{\rm max}\left\{|\sigma_{1}|,|\sigma_{2}|\right\}$,
 then $\supp \mathscr{F}u_{j}\subset D_{1}\cup D_{2}$ with $j=1,2$.
In this case, by using Lemma 2.5, since $-j+\frac{1}{2}\leq s\leq-\frac{j}{2},$ we have that
\begin{eqnarray*}
&&\left\|\Lambda^{-1}\partial_{x}(\prod_{j=1}^{2}u_{j})\right\|_{X^{s}}\leq C\left\|\langle k\rangle ^{-\frac{s}{j}}\langle \sigma \rangle ^{\frac{s}{j}}
(\mathscr{F}u_{1}*\mathscr{F}u_{2})\right\|_{l_{k}^{2}L_{\tau}^{2}}\nonumber\\&&
\leq C\left\|(J^{s}u_{1})(J^{s}u_{2})\right\|_{L_{xt}^{2}}\nonumber\\&&\leq
C\|u_{1}\|_{X_{s,\frac{1}{2}}}\|u_{2}\|_{X_{s,\frac{1}{2(2j+1)}+\epsilon}}\nonumber\\&&
\leq C\|u_{1}\|_{X_{s,\frac{1}{2}}}\|u_{2}\|_{X_{s,\frac{1}{2j}}}\leq C\prod_{j=1}^{2}\|u_{j}\|_{Z^{s}}.
\end{eqnarray*}
If $|\sigma|\leq 4{\rm max}\left\{|\sigma_{1}|,|\sigma_{2}|\right\},$ then we have $|\sigma|\sim |\sigma_{1}|$ or $|\sigma|\sim |\sigma_{2}|.$

\noindent When $|\sigma|\sim |\sigma_{1}|$,  $\supp\left( \mathscr{F}u_{1}*\mathscr{F}u_{2}\right) \subset D_{3},$ since $-j+\frac{1}{2}\leq s\leq-\frac{j}{2},$
 then we have that
\begin{eqnarray*}
&&\left\|\Lambda^{-1}\partial_{x}(\prod_{j=1}^{2}u_{j})\right\|_{X^{s}}\leq C\left\|\langle k\rangle ^{-\frac{s}{j}}\langle \sigma \rangle^{\frac{s}{j}}\left(\mathscr{F}u_{1}*\mathscr{F}u_{2}\right) \right\|_{l_{k}^{2}L_{\tau}^{2}}\nonumber\\
&&\leq C\left\|\left[\langle k\rangle^{-\frac{s}{j}-1}\langle\sigma\rangle^{\frac{s}{j}+1}\mathscr{F}u_{1}\right]*\left[\langle k\rangle^{-2j}\mathscr{F}u_{2}\right]\right\|_{l_{k}^{2}L_{\tau}^{2}}
\nonumber\\
&&\leq C\|u\|_{X_{-\frac{s}{j}-1,\frac{s}{j}+1}}\|\langle k\rangle^{-2j}\mathscr{F}u_{2}\|_{l_{k}^{1}L_{\tau}^{1}}\nonumber\\
&&\leq C\|u\|_{X_{-\frac{s}{j}-1,\frac{s}{j}+1}}\|\langle k\rangle^{s}\mathscr{F}u_{2}\|_{l_{k}^{2}L_{\tau}^{1}}
\leq
C\prod_{j=1}^{2}\|u_{j}\|_{Z^{s}}.
\end{eqnarray*}
When  $|\sigma|\sim |\sigma_{2}|$, this case can be proved similarly to case  $|\sigma|\sim |\sigma_{1}|$.

\noindent When (b)  occurs:  if  $|\sigma_{1}|>4{\rm max}\left\{|\sigma|,|\sigma_{2}|\right\}$
which yields $\supp\mathscr{F}u_{1}\subset D_{3} $   and    in this case, $\supp\left[\mathscr{F}u_{1}*\mathscr{F}u_{2}\right]\subset D_{1}\cup D_{2}$ and
$\supp\mathscr{F}u_{2}\subset D_{1}\cup D_{2}$, by using Lemmas 2.5,  2.3, since $-j+\frac{1}{2}\leq s\leq-\frac{j}{2},$ we have that
\begin{eqnarray*}
&&\left\|\Lambda^{-1}\partial_{x}(\prod_{j=1}^{2}u_{j})\right\|_{X^{s}}\leq C\left\|\langle k\rangle ^{s+1}\langle \sigma \rangle ^{-\frac{1}{2j}}\left(\mathscr{F}u_{1}*\mathscr{F}u_{2}\right)\right\|_{l_{k}^{2}L_{\tau}^{2}}\nonumber\\
&&\leq C\left\|(J^{-\frac{s}{j}-1}\Lambda ^{\frac{s}{j}+1}u_{1})(J^{-s-(2j-1)}u_{2})\right\|_{X_{0,-\frac{1}{2j}}}\leq C\|u_{1}\|_{X_{-\frac{s}{j}-1,\frac{s}{j}+1}}\|u_{2}\|_{X_{s,\frac{1}{2}}}\nonumber\\&&\leq
C\prod_{j=1}^{2}\|u_{j}\|_{Z^{s}}.
\end{eqnarray*}
If $|\sigma_{1}|\leq 4{\rm max}\left\{|\sigma|,|\sigma_{2}|\right\},$ we have $|\sigma_{1}|\sim |\sigma|$ or $|\sigma_{1}|\sim |\sigma_{2}|.$

\noindent
When $|\sigma_{1}|\sim |\sigma|$, this case can be proved similarly to $|\sigma|={\rm max}\left\{|\sigma|,|\sigma_{1}|,|\sigma_{2}|\right\}.$

\noindent
If $|\sigma_{1}|\sim |\sigma_{2}|,$   we can assume that $|\sigma|\leq C|k|^{2j+1},$
since $-j+\frac{1}{2}\leq s\leq -\frac{j}{2},$ we have
\begin{eqnarray*}
&&\left\|\Lambda^{-1}\partial_{x}(\prod_{j=1}^{2}u_{j})\right\|_{X^{s}}\leq C\left\|\langle k\rangle ^{s+1}\langle \sigma\rangle ^{-\frac{1}{2j}}
\left[\mathscr{F}u_{1}*\mathscr{F}u_{2}\right]\right\|_{l_{k}^{2}L_{\tau}^{2}}\nonumber\\
&&\leq C\left\|\langle k\rangle ^{s+\frac{3}{2}+\epsilon}\langle \sigma \rangle^{-\frac{1}{2j}+\frac{1}{2}+\epsilon}
\left[\mathscr{F}u_{1}*\mathscr{F}u_{2}\right]
\right\|_{l_{k}^{\infty}L_{\tau}^{\infty}}\nonumber\\
&&\leq C\left\|\left(\langle k\rangle ^{s+j+1-\frac{1}{2j}+(2j+2)\epsilon}\mathscr{F}u_{1}
\right)*\mathscr{F}u_{2}\right\|_{l_{k}^{\infty}l_{\tau}^{\infty}}\nonumber\\&&\leq C
\left\|\langle k\rangle ^{-3s+1-\frac{1}{2j}-3j+(2j+2)\epsilon}\right\|_{l_{k}^{\infty}}
\prod_{j=1}^{2}\|u_{j}\|_{X_{-\frac{s}{j}-1,\frac{s}{j}+1}}\nonumber\\&&\leq C
\prod_{j=1}^{2}\|u_{j}\|_{X_{-\frac{s}{j}-1,\frac{s}{j}+1}}\leq C\prod_{j=1}^{2}\|u_{j}\|_{Z^{s}}.
\end{eqnarray*}

We have completed the proof of Lemma 3.1.

\begin{Lemma}\label{Lemma3.2}
Let $j\geq 2$ and $-j+\frac{1}{2}\leq s\leq-\frac{j}{2}$. Then
\begin{eqnarray}
      \left\|\Lambda^{-1}\partial_{x}(\prod_{j=1}^{2}u_{j})\right\|_{Z^{s}}
      \leq C\prod\limits_{j=1}^{2}\|u_{j}\|_{Z^{s}}.
        \label{3.02}
\end{eqnarray}
\end{Lemma}
{\bf Proof.}
 Obviously, $\left(\R\times\dot{Z}_{\lambda}\right)^{2}\subset \bigcup\limits_{j=1}^{8}\Omega_{j},$
where $\Omega_{j}(1\leq j\leq 8)$ are defined as in Lemma 3.1.

\noindent (1) In region $\Omega_{1}$.
By using the Lemma 2.5 and the H\"older inequality as well as the Cauchy-Schwartz inequality, we have that
\begin{eqnarray*}
&&\left\|\Lambda^{-1}\partial_{x}(\prod_{j=1}^{2}u_{j})\right\|_{Y^{s}}\leq C\left\|\Lambda^{-1}\partial_{x}(\prod_{j=1}^{2}u_{j})\right\|_{X_{s,\frac{2j-1}{2j}}}\nonumber\\
&&\leq C\left\||k|\langle \sigma\rangle^{-\frac{1}{2j}}\left(\mathscr{F}u_{1}*\mathscr{F}u_{2}\right)
\right\|_{l_{k}^{2}L_{\tau}^{2}}\nonumber\\
&&\leq C\||k|\|_{l_{k}^{2}}\left\|\mathscr{F}u_{1}*\mathscr{F}u_{2}\right\|_{l_{k}^{\infty}L_{\tau}^{2}}\nonumber\\
&&\leq C\|\mathscr{F}u_{1}\|_{l_{k}^{2}L_{\tau}^{2}}\|\mathscr{F}u_{2}\|_{l_{k}^{2}L_{\tau}^{1}}\nonumber\\
&&\leq C\|u_{1}\|_{X_{s,\frac{1}{2j}}}\|u_{2}\|_{Y^{s}}\leq C\prod_{j=1}^{2}\|u_{j}\|_{Z^{s}}.
\end{eqnarray*}
(2) In region $\Omega_{2}$.
In this case, we consider (a)-(c) of Lemma 2.7, respectively.

\noindent (a) Case $|\sigma|={\rm max}\left\{|\sigma|,|\sigma_{1}|,|\sigma_{2}|\right\},$  by using the Young inequality and Lemma 2.7, since $-j+\frac{1}{2}\leq s\leq-\frac{j}{2},$
  we have that
\begin{eqnarray*}
&&\left\|\langle k\rangle^{s+1}\langle\sigma\rangle^{-1}
\left[\mathscr{F}u_{1}*\mathscr{F}u_{2}\right]\right\|_{l_{k}^{2}L_{\tau}^{1}}\nonumber\\
&&\leq C\left\|(\langle k\rangle ^{-j}\mathscr{F}u_{1})*(\langle k\rangle ^{-j}\mathscr{F}u_{2})\right\|_{l_{k}^{\infty}L_{\tau}^{1}}\leq C\prod_{j=1}^{2}\|u_{j}\|_{Y^{s}}.
\end{eqnarray*}
(b) Case $|\sigma_{1}|={\rm max}\left\{|\sigma|,|\sigma_{1}|,|\sigma_{2}|\right\},$
 we consider the following cases:
\begin{eqnarray*}
(i): |\sigma_{1}|>4{\rm max}\left\{|\sigma|,|\sigma_{2}|\right\},
(ii):|\sigma_{1}|\leq4{\rm max}\left\{|\sigma|,|\sigma_{2}|\right\},
\end{eqnarray*}
respectively.

\noindent
When (i) occurs:
if $\supp \mathscr{F}u_{1}\subset D_{1}$ which yields that $1\leq|k|\leq C$, by using Lemmas 2.5, 2.7, 2.3,   since $-j+\frac{1}{2}\leq s\leq-\frac{j}{2},$
 we have that
\begin{eqnarray*}
&&\left\|\Lambda^{-1}\partial_{x}(u_{1}u_{2})\right\|_{Y^{s}}\leq C\left\|\Lambda^{-1}\partial_{x}(u_{1}u_{2})\right\|_{Z^{s}}\nonumber\\&&
\leq C\left\|\partial_{x}(u_{1}u_{2})\right\|_{X_{s,-\frac{1}{2j}}}\nonumber\\&&
\leq C\left\|\langle k\rangle ^{s}\langle \sigma \rangle ^{-\frac{1}{2j}}
(\mathscr{F}u_{1}*\mathscr{F}u_{2})\right\|_{l_{k}^{2}L_{\tau}^{2}}\nonumber\\&&
\leq C\left\|\sigma \rangle ^{-\frac{1}{2j}}(\langle k\rangle^{s}\langle\sigma\rangle^{\frac{2j-1}{2j}}\mathscr{F}u_{1})*(\langle k\rangle ^{-s-2j+1}\mathscr{F}u_{2})\right\|_{l_{k}^{2}L_{\tau}^{2}}\nonumber\\
&&\leq C\left\|\left(J^{s}\Lambda ^{\frac{2j-1}{2j}}u_{1}\right)
\left(J^{-s-2j+1}u_{2}\right)\right\|_{X_{0,-\frac{1}{2j}}}\nonumber\\
&&\leq C\|u_{1}\|_{X_{s,\frac{2j-1}{2j}}}\|u_{2}\|_{X_{s,\frac{1}{2}}}\leq
C\prod_{j=1}^{2}\|u_{j}\|_{Z^{s}};
\end{eqnarray*}
if $\supp \mathscr{F}u_{1}\subset D_{2},$  by using Lemmas 2.5, 2.7, 2.3, since $-j+\frac{1}{2}\leq s\leq-\frac{j}{2}\leq -1$,
 we have that
\begin{eqnarray*}
&&\left\|\Lambda^{-1}\partial_{x}(u_{1}u_{2})\right\|_{Y^{s}}\leq C\left\|\Lambda^{-1}\partial_{x}(u_{1}u_{2})\right\|_{Z^{s}}\nonumber\\&&
\leq C\left\|\partial_{x}(u_{1}u_{2})\right\|_{X_{s,-\frac{1}{2j}}}\nonumber\\&&\leq \left\|\langle k\rangle ^{s+1}\langle \sigma \rangle ^{-\frac{1}{2j}}
\left[\mathscr{F}u_{1}*\mathscr{F}u_{2}\right]\right\|_{l_{k}^{2}L_{\tau}^{2}}\nonumber\\&&
\leq C\left\|\left(J^{(1-2j)s-1}\Lambda ^{s+1}u_{1}\right)\left(J^{-s-2j+1}u_{2}\right)
\right\|_{X_{0,-\frac{1}{2j}}}\nonumber\\
&&\leq C\|u_{1}\|_{X_{(1-2j)s-1,s+1}}\|u_{2}\|_{X_{s,\frac{1}{2}}}\leq
C\prod_{j=1}^{2}\|u_{j}\|_{Z^{s}}.
\end{eqnarray*}
When (ii) occurs: we have $|\sigma_{1}|\sim |\sigma|$ or $|\sigma_{1}|\sim |\sigma_{2}|$.

\noindent When $|\sigma_{1}|\sim |\sigma|$  is valid, this case can be proved similarly to
$|\sigma|={\rm max}\left\{|\sigma|,|\sigma_{1}|,|\sigma_{2}|\right\}.$
When $|\sigma_{1}|\sim |\sigma_{2}|$, if  $\supp \mathscr{F}u_{1}\subset D_{1}$ which leads to that $1\leq |k|\leq C$,
by using Lemmas 2.5, 2.7, 2.3, since $-j+\frac{1}{2}\leq s \leq -\frac{j}{2},$
 we have that
\begin{eqnarray*}
&&\left\|\Lambda^{-1}\partial_{x}(u_{1}u_{2})\right\|_{Y^{s}}\leq C
\left\|\Lambda^{-1}\partial_{x}(u_{1}u_{2})\right\|_{Z^{s}}
\leq C\left\|\partial_{x}(u_{1}u_{2})\right\|_{X_{s,-\frac{1}{2j}}}\nonumber\\
&&\leq C
\left\|\langle k\rangle ^{s+1}\langle \sigma \rangle ^{-\frac{1}{2j}}
(\mathscr{F}u_{1}*\mathscr{F}u_{2})\right\|_{l_{k}^{2}L_{\tau}^{2}}\nonumber\\
&&
\leq C\left\|\langle \sigma \rangle ^{-\frac{1}{2j}}\left[(\langle k\rangle^{s}\langle\sigma\rangle^{\frac{2j-1}{2j}}
\mathscr{F}u_{1})*(\langle k\rangle ^{-s-2j+1}\mathscr{F}u_{2})\right]
\right\|_{l_{k}^{2}L_{\tau}^{2}}\nonumber\\
&&\leq C\left\|\left(J^{s}\Lambda ^{\frac{2j-1}{2j}}u_{1}\right)
\left(J^{-s-2j+1}u_{2}\right)\right\|_{X_{0,-\frac{1}{2j}}}\nonumber\\&&\leq C
\|u_{1}\|_{X_{s,\frac{2j-1}{2j}}}\|u_{2}\|_{X_{s,\frac{1}{2}}}\leq
C\prod_{j=1}^{2}\|u_{j}\|_{Z^{s}};
\end{eqnarray*}
if $\supp \mathscr{F}u_{1}\subset D_{2},$ we  can assume that $\supp\mathscr{F} u_{2}\subset D_{2}$,
since $-j+\frac{1}{2}\leq s\leq-\frac{j}{2}\leq -1$,
 we can assume that  $|\sigma|\leq C|k_{1}|^{2j+1},$
 by using the H\"older  inequality  and the  Young inequality,
since $\langle k_{1}\rangle ^{(2j-1)s+1}\langle \sigma_{1}\rangle ^{-s-1}\leq
C\langle k_{1}\rangle ^{-2s-2j}$,  since $-j+\frac{1}{2}\leq s \leq -\frac{j}{2},$  we  have that
\begin{eqnarray*}
&&\left\|\Lambda^{-1}\partial_{x}(u_{1}u_{2})\right\|_{Y^{s}}\leq C
\left\|\Lambda^{-1}\partial_{x}(u_{1}u_{2})\right\|_{Z^{s}}
\leq C\left\|\partial_{x}(u_{1}u_{2})\right\|_{X_{s,-\frac{1}{2j}}}
\nonumber\\&&\leq C\left\|\langle k\rangle ^{s+1}\langle \sigma\rangle ^{-\frac{1}{2j}}
\left[\mathscr{F}u_{1}*\mathscr{F}u_{2}\right]\right\|_{l_{k}^{2}L_{\tau}^{2}}\nonumber\\
&&\leq C\left\|\langle k\rangle ^{s+\frac{3}{2}+\epsilon}\langle
\sigma \rangle^{-\frac{1}{2j}+\frac{1}{2}+\epsilon}
\left[\mathscr{F}u_{1}*\mathscr{F}u_{2}\right]
\right\|_{l_{k}^{\infty}L_{\tau}^{\infty}}\nonumber\\
&&\leq C\left\|\left(\langle k\rangle ^{s+j+1-\frac{1}{2j}+(2j+2)\epsilon}\mathscr{F}u_{1}
\right)*\mathscr{F}u_{2}\right\|_{l_{k}^{\infty}l_{\tau}^{\infty}}\nonumber\\&&\leq C
\left\|\langle k\rangle ^{-3s+1-\frac{1}{2j}-3j+(2j+2)\epsilon}\right\|_{l_{k}^{\infty}}
\prod_{j=1}^{2}\|u_{j}\|_{X_{(1-2j)s-1,s+1}}\nonumber\\&&\leq C
\prod_{j=1}^{2}\|u_{j}\|_{X_{(1-2j)s-1,s+1}}\leq C\prod_{j=1}^{2}\|u_{j}\|_{Z^{s}}.
\end{eqnarray*}
(c) Case $|\sigma_{2}|={\rm max}\left\{|\sigma|,|\sigma_{1}|,|\sigma_{2}|\right\}.$
 This case can be proved similarly to case (b).

\noindent (3) Region $\Omega_{3}$.
We  consider $|k|\leq |k_{1}|^{-2j}$  and  $|k_{1}|^{-2j}<|k|\leq 1,$ respectively.

\noindent When $|k|\leq |k_{1}|^{-2j}$, by using the Young inequality,  since $-j+\frac{1}{2}\leq s\leq-\frac{j}{2},$ we have that
\begin{eqnarray*}
&&\left\|\Lambda^{-1}\partial_{x}(\prod_{j=1}^{2}u_{j})\right\|_{Y^{s}}\leq C\left\||k|\langle \sigma\rangle ^{-\frac{1}{2j}}
\left[\mathscr{F}u_{1}*\mathscr{F}u_{2}\right]\right\|_{l_{k}^{2}L_{\tau}^{2}}\nonumber\\
&&\leq C\left\|\left[(\langle k\rangle^{-(2j-1)}\mathscr{F}u_{1})*(\langle k\rangle^{-(2j-1)}\mathscr{F}u_{2})\right]
\right\|_{l_{k}^{\infty}L_{\tau}^{2}}\nonumber\\
&&\leq C\|u_{1}\|_{X_{1-2j,0}}\|u_{2}\|_{Y^{1-2j}}\leq C\|u_{1}\|_{X_{1-2j,0}}\|u_{2}\|_{Y^{s}}
\leq
C\prod_{j=1}^{2}\|u_{j}\|_{Z^{s}}.
\end{eqnarray*}
When  $|k_{1}|^{-2j}<|k|\leq 1,$ we consider cases (a)-(c) of Lemma 2.7, respectively.

\noindent When (a) occurs:  by using the H\"older inequality and the Young inequality, since $-j+\frac{1}{2}\leq s\leq-\frac{j}{2},$ we have that
\begin{eqnarray*}
&&\left\|\Lambda^{-1}\partial_{x}(\prod_{j=1}^{2}u_{j})\right\|_{Y^{s}}\leq
 C\left\||k|\langle k\rangle^{s}\langle \sigma \rangle^{-1}\left[\mathscr{F}u_{1}*\mathscr{F}u_{2}\right]
\right\|_{l_{k}^{2}L_{\tau}^{1}}\nonumber\\
&&\leq C\left\|\langle k\rangle^{s}\left[\langle k\rangle ^{-j}\mathscr{F}u_{1}*\langle k\rangle ^{-j}\mathscr{F}u_{2}\right]
\right\|_{l_{k}^{2}L_{\tau}^{1}}\nonumber\\
&&\leq C\left\|\left[(\langle k\rangle ^{-j}\mathscr{F}u_{1})*(\langle k\rangle ^{-j}\mathscr{F}u_{2})\right]
\right\|_{l_{k}^{\infty}L_{\tau}^{1}}\nonumber\\
&&\leq C\prod_{j=1}^{2}\|\langle k\rangle^{-j}\mathscr{F}u_{j}\|_{l_{k}^{2}L_{\tau}^{1}}\nonumber\\
&&\leq C\prod_{j=1}^{2}\|u_{j}\|_{Y^{s}}\leq C\prod_{j=1}^{2}\|u_{j}\|_{Z^{s}}.
\end{eqnarray*}
When  (b) occurs:  we consider  $|\sigma_{1}|>4{\rm max}\left\{|\sigma|,|\sigma_{2}|\right\}$  and
$|\sigma_{1}|\leq4{\rm max}\left\{|\sigma|,|\sigma_{2}|\right\}$, respectively.

\noindent When $|\sigma_{1}|>4{\rm max}\left\{|\sigma|,|\sigma_{2}|\right\}$,
 we have $\supp \mathscr{F}u_{1}\subset D_{1}$
 and $X_{s,\frac{1}{2}+\epsilon}\hookrightarrow Y^{s},$
by using the H\"older inequality and the Young inequality,
since $-j+\frac{1}{2}\leq s\leq-\frac{j}{2},$ we have that
\begin{eqnarray*}
&&\left\|\Lambda^{-1}\partial_{x}(\prod_{j=1}^{2}u_{j})\right\|_{Y^{s}}\leq C\left\||k|\langle \sigma \rangle ^{-\frac{1}{2}+\epsilon}\left[\mathscr{F}u_{1}*\mathscr{F}u_{2}\right]\right\|_{l_{k}^{2}L_{\tau}^{2}}\nonumber\\
&&\leq C\left\||k|^{\frac{1}{2j}}\left(\langle k\rangle^{s}\langle \sigma \rangle ^{\frac{2j-1}{2j}}\mathscr{F}u_{1}\right)*\left(\langle k\rangle ^{-s-2j+1}\mathscr{F}u_{2}\right)\right\|_{l_{k}^{2}L_{\tau}^{2}}\nonumber\\
&&\leq C\left\|\left(\langle k\rangle^{s}\langle \sigma \rangle^{\frac{2j-1}{2j}}\mathscr{F}u_{1}\right)*\left(\langle k\rangle ^{-s-2j+1}\mathscr{F}u_{2}\right)\right\|_{l_{k}^{\infty}L_{\tau}^{2}}\leq C\|u_{1}\|_{X_{s,\frac{2j-1}{2j}}}\|u_{2}\|_{Y^{s}}.
\end{eqnarray*}
When $|\sigma_{1}|\leq4{\rm max}\left\{|\sigma|,|\sigma_{2}|\right\}$, we have that $|\sigma_{1}|\sim |\sigma|$ or $|\sigma_{1}|\sim |\sigma_{2}|$.

\noindent When $|\sigma_{1}|\sim |\sigma|$, this case can be proved similarly to case $|\sigma|={\rm max}\left\{|\sigma|,|\sigma_{1}|,|\sigma_{2}|\right\}.$
When $|\sigma_{1}|\sim |\sigma_{2}|$, if $\supp \mathscr{F}u_{j}\subset D_{1}$ with $j=1,2,$ by using $X_{s,\frac{1}{2}+\epsilon}\hookrightarrow Y^{s}$ and the Young inequality, since $-j+\frac{1}{2}\leq s \leq -\frac{j}{2},$  we have that
\begin{eqnarray*}
&&\left\|\Lambda^{-1}\partial_{x}(\prod_{j=1}^{2}u_{j})\right\|_{Y^{s}}\leq C\left\||k|\langle \sigma \rangle ^{-\frac{1}{2}+\epsilon}\left[\mathscr{F}u_{1}*\mathscr{F}u_{2}\right]
\right\|_{l_{k}^{2}L_{\tau}^{2}}\nonumber\\
&&\leq C\left\||k|^{\frac{1}{2j}}\langle \sigma \rangle ^{-\frac{1}{2}+\epsilon}\left(\langle k\rangle^{s}\langle \sigma \rangle ^{\frac{2j-1}{2j}}\mathscr{F}u_{1}\right)*\left(\langle k\rangle ^{-s-2j+1}\mathscr{F}u_{2}\right)\right\|_{l_{k}^{2}L_{\tau}^{2}}\nonumber\\
&&\leq C\left\|\langle \sigma \rangle ^{-\frac{1}{2}+\epsilon}\left(\langle k\rangle^{s}\langle \sigma \rangle^{\frac{2j-1}{2j}}
\mathscr{F}u_{1}\right)*\left(\langle k\rangle ^{-s-2j+1}\mathscr{F}u_{2}\right)
\right\|_{l_{k}^{2}L_{\tau}^{2}}\nonumber\\
&&\leq C\|u_{1}\|_{X_{s,\frac{2j-1}{2j}}}\|u_{2}\|_{X_{s,\frac{1}{2j}}}\leq C\prod_{j=1}^{2}\|u_{j}\|_{Z^{s}};
\end{eqnarray*}
if $\supp \mathscr{F}u_{j}\subset D_{2},$ by using $X_{s,\frac{1}{2}+\epsilon}\hookrightarrow Y^{s}$ and the Young inequality, since $-j+\frac{1}{2}\leq s \leq -\frac{j}{2},$  we have that
\begin{eqnarray*}
&&\left\|\Lambda^{-1}\partial_{x}(\prod_{j=1}^{2}u_{j})\right\|_{Y^{s}}\leq
 C\left\||k|\langle \sigma \rangle ^{-\frac{1}{2}+\epsilon}\left[\mathscr{F}u_{1}*\mathscr{F}u_{2}\right]
 \right\|_{l_{k}^{2}L_{\tau}^{2}}\nonumber\\
&&\leq C\left\|\left(\langle k\rangle^{(1-2j)s-1}\langle \sigma \rangle ^{s+1}
\mathscr{F}u_{1}\right)*\left(\langle k\rangle ^{-2s-2j}\mathscr{F}u_{2}\right)\right\|_{X^{0,-\frac{1}{2}+\epsilon}}\nonumber\\
&&\leq C\left\|u_{1}\right\|_{X^{(1-2j)s,s+1}}\|u\|_{X_{-2s-2j,\frac{s}{j}+1}}\nonumber\\
&&\leq C\|u_{1}\|_{X_{(1-2j)s,s+1}}\|u_{2}\|_{X_{-2s-2j,\frac{s}{j}+1}}\nonumber\\
&&\leq C\prod_{j=1}^{2}\|u_{j}\|_{X_{(1-2j)s,s+1}}\leq C\prod_{j=1}^{2}\|u_{j}\|_{Z^{s}};
\end{eqnarray*}
if $\supp \mathscr{F}u_{j}\subset D_{3},$ by using $X_{s,\frac{1}{2}+\epsilon}\hookrightarrow Y^{s}$ and the Young inequality, since $-j+\frac{1}{2}\leq s \leq -\frac{j}{2},$  we have that
\begin{eqnarray*}
&&\left\|\Lambda^{-1}\partial_{x}(\prod_{j=1}^{2}u_{j})\right\|_{Y^{s}}\leq
 C\left\||k|\langle \sigma \rangle ^{-\frac{1}{2}+\epsilon}\left[\mathscr{F}u_{1}*\mathscr{F}u_{2}\right]
 \right\|_{l_{k}^{2}L_{\tau}^{2}}\nonumber\\
&&\leq C\left\|\left(\langle k\rangle^{-\frac{s}{j}-1}\langle \sigma \rangle ^{\frac{s}{j}+1}
\mathscr{F}u_{1}\right)*\left(\langle k\rangle ^{-2s-2j}\mathscr{F}u_{2}\right)\right\|_{X^{0,-\frac{1}{2}+\epsilon}}\nonumber\\
&&\leq C\left\|u_{1}\right\|_{X^{-\frac{s}{j}-1,\frac{s}{j}+1}}\|u\|_{X_{-2s-2j,\frac{s}{j}+1}}\nonumber\\
&&\leq C\|u_{1}\|_{X_{-\frac{s}{j}-1,-\frac{s}{j}-1}}\|u_{2}\|_{X_{-2s-2j,\frac{s}{j}+1}}\nonumber\\
&&\leq C\prod_{j=1}^{2}\|u_{j}\|_{X_{-\frac{s}{j}-1,\frac{s}{j}+1}}\leq C\prod_{j=1}^{2}\|u_{j}\|_{Z^{s}}.
\end{eqnarray*}
(4) Region $\Omega_{4}$.
We consider  cases (a)-(c) of Lemma 2.7, respectively.

\noindent When (a) occurs: we consider case $|\sigma|>4{\rm max}\left\{|\sigma_{1}|,|\sigma_{2}|\right\}$ and $|\sigma|\leq 4{\rm max}\left\{|\sigma_{1}|,|\sigma_{2}|\right\},$ respectively.

\noindent If $|\sigma|>4{\rm max}\left\{|\sigma_{1}|,|\sigma_{2}|\right\}$,
 then $\supp \mathscr{F}u_{j}\subset D_{1}\cup D_{2}$ with $j=1,2$  and  $\supp \left[\mathscr{F}u_{1}*\mathscr{F}u_{2}\right]\subset D_{2}.$
In this case, by using Lemma 2.5, since $-j+\frac{1}{2}\leq s\leq-\frac{j}{2},$ we have that
\begin{eqnarray*}
&&\left\|\Lambda^{-1}\partial_{x}(u_{1}u_{2})\right\|_{Y^{s}}\nonumber\\&&\leq C\left\|\langle k\rangle^{s+1}\langle \sigma \rangle^{-1}
\left[\mathscr{F}u_{1}*\mathscr{F}u_{2}\right]\right\|_{l_{k}^{2}L_{\tau}^{1}}
 \nonumber\\
 &&\leq C\left\|\langle k\rangle^{s+1-2j}|k_{1}|^{-1}
\left[\mathscr{F}u_{1}*\mathscr{F}u_{2}\right]\right\|_{l_{k}^{2}L_{\tau}^{1}}
 \nonumber\\
 &&\leq C\left\|
\left[(\langle k\rangle^{-2j}\mathscr{F}u_{1})*(\langle k\rangle^{s}\mathscr{F}u_{2})\right]\right\|_{l_{k}^{2}L_{\tau}^{1}}
 \nonumber\\
&&\leq C\|\langle k\rangle^{-2j}\mathscr{F}u_{1}\|_{l_{k}^{1}L_{\tau}^{1}}\|u_{2}\|_{l_{k}^{2}L_{\tau}^{1}}
\nonumber\\&&\leq C\prod_{j=1}^{2}\|u_{j}\|_{Y^{s}}
\leq C\prod_{j=1}^{2}\|u_{j}\|_{Z^{s}}.
\end{eqnarray*}
When $|\sigma|\leq4{\rm max}\left\{|\sigma_{1}|,|\sigma_{2}|\right\}$, we have  $|\sigma|\sim |\sigma_{1}|$ or $|\sigma|\sim |\sigma_{2}|$.

\noindent
 When  $|\sigma|\sim |\sigma_{1}|$,
if $\supp \left[\mathscr{F}u_{1}*\mathscr{F}u_{2}\right]\subset D_{2},$
 by using the Young inequality,  we have that
\begin{eqnarray*}
&&\left\|\langle k\rangle^{s+1}\langle \sigma \rangle^{-1}
\left[\mathscr{F}u_{1}*\mathscr{F}u_{2}\right]\right\|_{l_{k}^{2}L_{\tau}^{1}} \nonumber\\&&\leq\left\|(\langle k\rangle^{-2j}\mathscr{F}u_{1})*(\langle k\rangle^{s}\mathscr{F}u_{2})\right\|_{L_{k}^{2}l_{\tau}^{1}}\nonumber\\&&\leq C\|\langle k\rangle^{-2j}\mathscr{F}u_{1}\|_{l_{k}^{1}L_{\tau}^{2}}\|u_{2}\|_{Y_{s}}\leq C\prod_{j=1}^{2}\|u_{j}\|_{Z^{s}};
\end{eqnarray*}
if $\supp \left[\mathscr{F}u_{1}*\mathscr{F}u_{2}\right]\subset D_{3},$
 $\supp \mathscr{F}u_{1}\subset D_{3},$ by using Lemma 2.3, since $-j+\frac{1}{2}\leq s\leq-\frac{j}{2},$ we have that
\begin{eqnarray*}
&&\left\|\Lambda^{-1}\partial_{x}(u_{1}u_{2})\right\|_{Y^{s}}\leq C\left\|\Lambda^{-1}\partial_{x}(u_{1}u_{2})\right\|_{Z^{s}}
\leq C\left\|\partial_{x}(u_{1}u_{2})\right\|_{X_{s,-\frac{1}{2j}}}\nonumber\\&&\leq C\left\|\langle k\rangle^{s+1}\langle \sigma \rangle^{-\frac{1}{2j}}
\left[\mathscr{F}u_{1}*\mathscr{F}u_{2}\right]\right\|_{l_{k}^{2}L_{\tau}^{2}} \nonumber\\&&\leq\left\|(J^{-\frac{s}{j}-1}\Lambda^{\frac{s}{j}+1}u_{1})(J^{-s-(2j-1)}u_{2})\right\|_{X_{0,-\frac{1}{2j}}}\nonumber\\&&\leq C\|u_{1}\|_{X_{-\frac{s}{j}-1,\frac{s}{j}+1}}\|u_{2}\|_{X_{s,\frac{1}{2}}}\leq C\prod_{j=1}^{2}\|u_{j}\|_{Z^{s}}.
\end{eqnarray*}
When $|\sigma|\sim |\sigma_{2}|$, this case can be proved similarly to case   $|\sigma|\sim |\sigma_{1}|$.

\noindent (b): $|\sigma_{1}|={\rm max}\left\{|\sigma|,|\sigma_{1}|,|\sigma_{2}|\right\}.$
If $|\sigma_{1}|>4{\rm max}\left\{|\sigma|,|\sigma_{2}|\right\}$,
 then $\supp \mathscr{F}u_{1}\subset D_{3}$ and  $\supp \mathscr{F}u_{2}\subset D_{1}\cup D_{2}$,
In this case, by using Lemma 2.3,   since $-j+\frac{1}{2}\leq s\leq-\frac{j}{2},$  we have that
\begin{eqnarray*}
&&\left\|\Lambda^{-1}\partial_{x}(u_{1}u_{2})\right\|_{Y^{s}}\leq C\left\|\Lambda^{-1}\partial_{x}(u_{1}u_{2})\right\|_{Z^{s}}
\leq C\left\|\partial_{x}(u_{1}u_{2})\right\|_{X_{s,-\frac{1}{2j}}}\nonumber\\&&\leq C\left\|\langle k\rangle^{s+1}\langle \sigma \rangle^{-\frac{1}{2j}}
\left[\mathscr{F}u_{1}*\mathscr{F}u_{2}\right]\right\|_{l_{k}^{2}L_{\tau}^{2}}
\nonumber\\
&&\leq C\left\|(J^{-\frac{s}{j}-1}\langle \sigma \rangle ^{\frac{s}{j}+1}u_{1})(J^{s-(2j-1)}u_{2})\right\|_{X_{0,-\frac{1}{2j}}}
\nonumber\\&&\leq C\|u_{1}\|_{X_{-\frac{s}{j}-1,\frac{s}{j}+1}}\|u_{2}\|_{X_{s,\frac{1}{2}}}
\leq C\prod_{j=1}^{2}\|u_{j}\|_{Z^{s}}.
\end{eqnarray*}
When $|\sigma_{1}|\leq4{\rm max}\left\{|\sigma|,|\sigma_{2}|\right\}$,
we have that $|\sigma_{1}|\sim |\sigma|$ or $|\sigma_{1}|\sim |\sigma_{2}|$.

\noindent
 When  $|\sigma_{1}|\sim |\sigma|$, this case can be proved similarly
 to case $|\sigma|={\rm max}\left\{|\sigma|,|\sigma_{1}|,|\sigma_{2}|\right\}.$

\noindent When $|\sigma_{1}|\sim |\sigma_{2}|$, we have $\supp \mathscr{F}u_{1}\subset D_{3}$
and  $\supp \mathscr{F}u_{2}\subset D_{2}\bigcup D_{3}.$

\noindent
When  $\supp \mathscr{F}u_{2}\subset D_{2}$, by using Lemmas 2.5,  2.3, since $-j+\frac{1}{2}\leq s\leq -\frac{j}{2},$  we have that
\begin{eqnarray*}
&&\left\|\Lambda^{-1}\partial_{x}(u_{1}u_{2})\right\|_{Y^{s}}\leq C
\left\|\Lambda^{-1}\partial_{x}(u_{1}u_{2})\right\|_{Z^{s}}
\leq C\left\|\partial_{x}(u_{1}u_{2})\right\|_{X_{s,-\frac{1}{2j}}}\nonumber\\&&\leq C\left\|\langle k\rangle^{s}\langle \sigma \rangle^{-\frac{1}{2j}}
\left[(|k|\mathscr{F}u_{1})*\mathscr{F}u_{2}\right]\right\|_{l_{k}^{2}L_{\tau}^{2}} \nonumber\\&&\leq C\left\|\langle k\rangle ^{s+\frac{1}{2}+\epsilon}\langle \sigma \rangle^{-\frac{1}{2j}+\frac{1}{2}+\epsilon}
\left[(|k|\mathscr{F}u_{1})*\mathscr{F}u_{2}\right]
\right\|_{l_{k}^{\infty}L_{\tau}^{\infty}}\nonumber\\
&&\leq C\left\|\langle k\rangle ^{s+j-\frac{1}{2j}+(2j+2)\epsilon}\left((|k|\mathscr{F}u_{1})
*\mathscr{F}u_{2}\right)\right\|_{l_{k}^{\infty}l_{\tau}^{\infty}}\nonumber\\&&\leq C
\left\|\langle k\rangle ^{-3s+1-\frac{1}{2j}-3j+(2j+2)\epsilon}\right\|_{l_{k}^{\infty}}
\|u_{1}\|_{X_{-\frac{s}{j}-1,\frac{s}{j}+1}}\|u_{2}\|_{X_{(1-2j)s-1,s+1}}\nonumber\\&&\leq C
\prod_{j=1}^{2}\|u_{j}\|_{Z^{s}}.
\end{eqnarray*}
When   $\supp \mathscr{F}u_{2}\subset D_{3}$,
 by using Lemma 2.5, since $-j+\frac{1}{2}\leq s\leq-\frac{j}{2},$ we have that
\begin{eqnarray*}
&&\left\|\Lambda^{-1}\partial_{x}(u_{1}u_{2})\right\|_{Y^{s}}\leq C
\left\|\Lambda^{-1}\partial_{x}(u_{1}u_{2})\right\|_{Z^{s}}
\leq C\left\|\partial_{x}(u_{1}u_{2})\right\|_{X_{s,-\frac{1}{2j}}}\nonumber\\&&\leq C\left\|\langle k\rangle ^{s+1}\langle \sigma\rangle ^{-\frac{1}{2j}}
\left[\mathscr{F}u_{1}*\mathscr{F}u_{2}\right]\right\|_{l_{k}^{2}L_{\tau}^{2}}\nonumber\\
&&\leq C\left\|\langle k\rangle ^{s+\frac{3}{2}+\epsilon}\langle \sigma \rangle^{-\frac{1}{2j}+\frac{1}{2}+\epsilon}
\left[\mathscr{F}u_{1}*\mathscr{F}u_{2}\right]
\right\|_{l_{k}^{\infty}L_{\tau}^{\infty}}\nonumber\\
&&\leq C\left\|\left(\langle k\rangle ^{s+j+1-\frac{1}{2j}+(2j+2)\epsilon}\mathscr{F}u_{1}
\right)*\mathscr{F}u_{2}\right\|_{l_{k}^{\infty}l_{\tau}^{\infty}}\nonumber\\&&\leq C
\left\|\langle k\rangle ^{-3s+1-\frac{1}{2j}-3j+(2j+2)\epsilon}\right\|_{l_{k}^{\infty}}
\prod_{j=1}^{2}\|u_{j}\|_{X_{-\frac{s}{j}-1,\frac{s}{j}+1}}\nonumber\\&&\leq C
\prod_{j=1}^{2}\|u_{j}\|_{X_{-\frac{s}{j}-1,\frac{s}{j}+1}}\leq C\prod_{j=1}^{2}\|u_{j}\|_{Z^{s}}.
\end{eqnarray*}
(5) In region $\Omega_{5}$.
 In this region, we consider the case $|k_{1}|\leq |k|^{-2j}$ and $|k|^{-2j}< |k_{1}|\leq 1,$ respectively.

\noindent When $|k_{1}|\leq |k|^{-2j}$, by using Lemma 2.5, the Young inequality and Cauchy-Schwartz inequality, since $-j+\frac{1}{2}\leq s\leq-\frac{j}{2},$ we have that
\begin{eqnarray*}
&&\left\|\Lambda^{-1}\partial_{x}(\prod_{j=1}^{2}u_{j})\right\|_{Y^{s}}\leq C\left\|\langle k\rangle^{s+1}\langle \sigma\rangle ^{-\frac{1}{2j}}
\left[\mathscr{F}u_{1}*\mathscr{F}u_{2}\right]\right\|_{l_{k}^{2}L_{\tau}^{2}}\nonumber\\
&&\leq C\left\|\left[(|k|^{-\frac{1}{2j}}\mathscr{F}u_{1})*(\langle k\rangle^{s}\mathscr{F}u_{2})\right]\right\|_{l_{k}^{2}L_{\tau}^{2}}\nonumber\\
&&\leq C\left\||k|^{-\frac{1}{2j}}\mathscr{F}u_{1}\right\|_{l_{k}^{1}l_{\tau}^{2}}\|u_{2}\|_{Y^{s}}\nonumber\\&&
\leq C\|\mathscr{F}u_{1}\|_{l_{n}^{2}L_{\tau}^{2}}\|u_{2}\|_{Y^{s}}
\leq C\prod_{j=1}^{2}\|u_{j}\|_{Z^{s}}.
\end{eqnarray*}
When  $|k|^{-2j}< |k_{1}|\leq 1.$
We consider   cases (a)-(c) of Lemma 2.7, respectively.

\noindent When (a) occurs:
 by using Lemma 2.5, the Young inequality and Cauchy-Schwartz inequality, since $-j+\frac{1}{2}\leq s\leq-\frac{j}{2},$ we have that
\begin{eqnarray*}
&&\left\|\Lambda^{-1}\partial_{x}(\prod_{j=1}^{2}u_{j})\right\|_{Y^{s}}\leq C\left\|\langle k\rangle^{s+1}\langle \sigma\rangle ^{-\frac{1}{2j}}
\left[\mathscr{F}u_{1}*\mathscr{F}u_{2}\right]\right\|_{l_{k}^{2}L_{\tau}^{2}}\nonumber\\
&&\leq C\left\|\left[(|k|^{-\frac{1}{2j}}\mathscr{F}u_{1})*(\langle k\rangle ^{s}\mathscr{F}u_{2})\right]
\right\|_{l_{k}^{2}L_{\tau}^{2}}\nonumber\\
&&\leq C\left\||k|^{-\frac{1}{2j}}\mathscr{F}u_{1}\right\|_{l_{k}^{1}l_{\tau}^{2}}\|u_{2}\|_{Y^{s}}\nonumber\\&&
\leq C\|\mathscr{F}u_{1}\|_{l_{k}^{2}L_{\tau}^{2}}\|u_{2}\|_{Y^{s}}
\leq C\prod_{j=1}^{2}\|u_{j}\|_{Z^{s}}.
\end{eqnarray*}
When (b) occurs:
 by using the Lemma 2.5, Young inequality and Cauchy-Schwartz inequality, since $-j+\frac{1}{2}\leq s\leq-\frac{j}{2},$ we have that
\begin{eqnarray*}
&&\left\|\Lambda^{-1}\partial_{x}(\prod_{j=1}^{2}u_{j})\right\|_{Y^{s}}\leq C\left\|\langle k\rangle ^{s+1}\langle \sigma\rangle ^{-\frac{1}{2j}}
\left[\mathscr{F}u_{1}*\mathscr{F}u_{2}\right]\right\|_{l_{k}^{2}L_{\tau}^{2}}\nonumber\\
&&\leq C\left\|\left[(|k|^{-\frac{1}{2j}}\langle \sigma \rangle ^{\frac{1}{2j}}\mathscr{F}u_{1})*(\langle k\rangle ^{s}\mathscr{F}u_{2})\right]
\right\|_{l_{k}^{2}L_{\tau}^{2}}\nonumber\\
&&\leq C\left\||k|^{-\frac{1}{2j}}\langle \sigma \rangle ^{\frac{1}{2j}}\mathscr{F}u_{1}\right\|_{l_{k}^{1}l_{\tau}^{2}}\|u_{2}\|_{Y^{s}}\nonumber\\&&
\leq C\|\langle \sigma \rangle ^{\frac{1}{2j}}\mathscr{F}u_{1}\|_{l_{k}^{2}L_{\tau}^{2}}\|u_{2}\|_{Y^{s}}
\leq C\|u_{1}\|_{X_{0,\frac{1}{2j}}}\|u_{2}\|_{Y^{s}}\nonumber\\&&
\leq C\prod_{j=1}^{2}\|u_{j}\|_{Z^{s}}.
\end{eqnarray*}
When (c)  occurs:
by using Lemma 2.5, the Young inequality and Cauchy-Schwartz inequality, since $-j+\frac{1}{2}\leq s\leq-\frac{j}{2},$ we have that
\begin{eqnarray*}
&&\left\|\Lambda^{-1}\partial_{x}(\prod_{j=1}^{2}u_{j})\right\|_{Y^{s}}\leq C\left\|(1-\partial_{x}^{2})^{-\frac{1}{2}}\partial_{x}(\prod_{j=1}^{2}u_{j})\right\|_{X^{s}}\leq C\left\|\langle k\rangle ^{s+1}\langle \sigma\rangle ^{-\frac{1}{2j}}
\left[\mathscr{F}u_{1}*\mathscr{F}u_{2}\right]\right\|_{l_{k}^{2}L_{\tau}^{2}}\nonumber\\
&&\leq C\left\|\left[(|k|^{-\frac{1}{2j}}\mathscr{F}u_{1})*(\langle k\rangle ^{s}\langle \sigma \rangle ^{\frac{1}{2j}}\mathscr{F}u_{2})\right]
\right\|_{l_{k}^{2}L_{\tau}^{2}}\nonumber\\
&&\leq C\left\||k|^{-\frac{1}{2j}}\mathscr{F}u_{1}\right\|_{l_{k}^{1}l_{\tau}^{1}}\|u_{2}\|_{X_{s,\frac{1}{2j}}}\nonumber\\&&
\leq C\|\mathscr{F}u_{1}\|_{l_{k}^{2}L_{\tau}^{1}}\|u_{2}\|_{X_{s,\frac{1}{2j}}}\nonumber\\&&
\leq C\|u_{1}\|_{Y^{s}}\|u_{2}\|_{X_{s,\frac{1}{2j}}}
\leq C\prod_{j=1}^{2}\|u_{j}\|_{Z^{s}}.
\end{eqnarray*}
(6)In region $\Omega_{6}$.
This region can be proved similarly to $\Omega_{4}$.

\noindent(7)In region $\Omega_{7}$.
This region can be proved similarly to $\Omega_{7}$.

\noindent (8)In region $\Omega_{8}$. We consider  cases (a)-(c) of Lemma 2.7, respectively.

\noindent
When (a) occurs: $\supp\left( \mathscr{F}u_{1}*\mathscr{F}u_{2}\right) \subset D_{3}$.
If $|\sigma|>4{\rm max}\left\{|\sigma_{1}|,|\sigma_{2}|\right\}$,
 then $\supp \mathscr{F}u_{j}\subset D_{1}\cup D_{2}$ with $j=1,2$.
In this case, by using Lemmas 2.5,  2.3, since $-j+\frac{1}{2}\leq s\leq-\frac{j}{2},$ we have that
\begin{eqnarray*}
&&\left\|\Lambda^{-1}\partial_{x}(u_{1}u_{2})\right\|_{Y^{s}}\nonumber\\&&\leq C\left\|\langle k\rangle ^{s+1}\langle \sigma \rangle ^{-1}
(\mathscr{F}u_{1}*\mathscr{F}u_{2})\right\|_{l_{k}^{2}L_{\tau}^{1}}\nonumber\\&&
\leq C\left\|\langle k\rangle ^{s-2j}
(\mathscr{F}u_{1}*\mathscr{F}u_{2})\right\|_{l_{k}^{2}L_{\tau}^{1}}\nonumber\\&&
\leq C\|\langle k\rangle ^{-2j}u_{1}\|_{L_{k}^{1}L_{\tau}^{1}}\|\langle k\rangle ^{s}u_{1}\|_{L_{k}^{2}L_{\tau}^{1}}\nonumber\\&&\leq C\prod_{j=1}^{2}\|u_{j}\|_{Y^{s}}
\leq C\prod_{j=1}^{2}\|u_{j}\|_{Z^{s}}.
\end{eqnarray*}
If $|\sigma|\leq 4{\rm max}\left\{|\sigma_{1}|,|\sigma_{2}|\right\},$  we have $|\sigma|\sim |\sigma_{1}|$ or $|\sigma|\sim |\sigma_{2}|.$

\noindent
When $|\sigma|\sim |\sigma_{1}|$. In this case, $\supp\left( \mathscr{F}u_{1}*\mathscr{F}u_{2}\right) \subset D_{3},$  by using $X_{s,\frac{1}{2}+\epsilon}\hookrightarrow Y^{s},$ since $-j+\frac{1}{2}\leq s\leq -\frac{j}{2},$
 then we have that
\begin{eqnarray*}
&&\left\|\Lambda^{-1}\partial_{x}(u_{1}u_{2})\right\|_{Y^{s}}
\leq C\left\|\langle k\rangle ^{s+1}\langle \sigma \rangle^{-1}\left(\mathscr{F}u_{1}*\mathscr{F}u_{2}\right) \right\|_{l_{k}^{2}L_{\tau}^{2}}\nonumber\\
&&\leq C\left\|(\langle k\rangle^{-2j}\mathscr{F}u_{1})*(\langle k\rangle^{s}\mathscr{F}u_{2})\right\|_{l_{k}^{2}L_{\tau}^{2}}\nonumber\\
&&\leq C\left\|\left[\langle k\rangle^{-\frac{s}{j}-1}\langle\sigma\rangle^{\frac{s}{j}+1}\mathscr{F}u_{1}\right]*\left[\langle k\rangle^{-s-4j}\mathscr{F}u_{2}\right]\right\|_{l_{k}^{2}L_{\tau}^{1}}
\nonumber\\
&&\leq C\|u_{1}\|_{X_{-\frac{s}{j}-1,\frac{s}{j}+1}}\|\langle k\rangle^{-s-4j}\mathscr{F}u_{2}\|_{l_{k}^{1}L_{\tau}^{1}}\nonumber\\
&&\leq C\|u_{1}\|_{X_{-\frac{s}{j}-1,\frac{s}{j}+1}}\|\langle k\rangle^{s}\mathscr{F}u_{2}\|_{l_{k}^{2}L_{\tau}^{1}}\leq
C\prod_{j=1}^{2}\|u_{j}\|_{Z^{s}}.
\end{eqnarray*}
When  $|\sigma|\sim |\sigma_{2}|$, this case can be proved similarly to case  $|\sigma|\sim |\sigma_{1}|$.

\noindent When (b)  occurs:  if  $|\sigma_{1}|>4{\rm max}\left\{|\sigma|,|\sigma_{2}|\right\}$
which yields $\supp\mathscr{F}u_{1}\subset D_{3} $   and    in this case, $\supp\left[\mathscr{F}u_{1}*\mathscr{F}u_{2}\right]\subset D_{1}\cup D_{2}$ and
$\mathscr{F}u_{2}\subset D_{1}\cup D_{2}$, by using Lemmas 2.5, 2.3, since $-j+\frac{1}{2}\leq s\leq-\frac{j}{2},$ we have that
\begin{eqnarray*}
&&\left\|\Lambda^{-1}\partial_{x}(u_{1}u_{2})\right\|_{Y^{s}}\leq C
\left\|\Lambda^{-1}\partial_{x}(u_{1}u_{2})\right\|_{Z^{s}}
\leq C\left\|\partial_{x}(u_{1}u_{2})\right\|_{X_{s,-\frac{1}{2j}}}\nonumber\\&&\leq C\left\|\langle k\rangle ^{s+1}\langle \sigma \rangle ^{-\frac{1}{2j}}\left(\mathscr{F}u_{1}*\mathscr{F}u_{2}\right)\right\|_{l_{k}^{2}L_{\tau}^{2}}\nonumber\\
&&\leq C\left\|(J^{-\frac{s}{j}-1}\Lambda ^{\frac{s}{j}+1}u_{1})(J^{-s-(2j-1)}u_{2})\right\|_{X_{0,-\frac{1}{2j}}}\leq C\|u_{1}\|_{X_{-\frac{s}{j}-1,\frac{s}{j}+1}}\|u_{2}\|_{X_{s,\frac{1}{2}}}\nonumber\\
&&\leq
C\prod_{j=1}^{2}\|u_{j}\|_{Z^{s}}.
\end{eqnarray*}
If $|\sigma_{1}|\leq 4{\rm max}\left\{|\sigma|,|\sigma_{2}|\right\},$ we have  $|\sigma_{1}|\sim |\sigma|$ or $|\sigma_{1}|\sim |\sigma_{2}|.$

\noindent
When $|\sigma_{1}|\sim |\sigma|$, this case can be proved similarly to $|\sigma|={\rm max}\left\{|\sigma|,|\sigma_{1}|,|\sigma_{2}|\right\}.$

\noindent
If $|\sigma_{1}|\sim |\sigma_{2}|,$  in this case, we can assume that $|\sigma|\leq C|k|^{2j+1},$
 since $-j+\frac{1}{2}\leq s\leq-\frac{j}{2}$ and
 $\epsilon<\frac{1}{100j},$  we have that
\begin{eqnarray*}
&&\left\|\Lambda^{-1}\partial_{x}(u_{1}u_{2})\right\|_{Y^{s}}\leq C
\left\|\Lambda^{-1}\partial_{x}(u_{1}u_{2})\right\|_{Z^{s}}
\leq C\left\|\partial_{x}(u_{1}u_{2})\right\|_{X_{s,-\frac{1}{2j}}}\nonumber\\&&\leq C\left\|\langle k\rangle ^{s+1}\langle \sigma\rangle ^{-\frac{1}{2j}}
\left[\mathscr{F}u_{1}*\mathscr{F}u_{2}\right]\right\|_{l_{k}^{2}L_{\tau}^{2}}\nonumber\\
&&\leq C\left\|\langle k\rangle ^{s+\frac{3}{2}+\epsilon}\langle \sigma \rangle^{-\frac{1}{2j}+\frac{1}{2}+\epsilon}
\left[\mathscr{F}u_{1}*\mathscr{F}u_{2}\right]
\right\|_{l_{k}^{\infty}L_{\tau}^{\infty}}\nonumber\\
&&\leq C\left\|\left(\langle k\rangle ^{s+j+1-\frac{1}{2j}+(2j+2)\epsilon}\mathscr{F}u_{1}
\right)*\mathscr{F}u_{2}\right\|_{l_{k}^{\infty}l_{\tau}^{\infty}}\nonumber\\&&\leq C
\left\|\langle k\rangle ^{-3s+1-\frac{1}{2j}-3j+(2j+2)\epsilon}\right\|_{l_{k}^{\infty}}
\prod_{j=1}^{2}\|u_{j}\|_{X_{-\frac{s}{j}-1,\frac{s}{j}+1}}\nonumber\\&&\leq C
\prod_{j=1}^{2}\|u_{j}\|_{X_{-\frac{s}{j}-1,\frac{s}{j}+1}}\leq C\prod_{j=1}^{2}\|u_{j}\|_{Z^{s}}.
\end{eqnarray*}

We have completed the proof of Lemma 3.2.

\begin{Lemma}\label{Lemma3.3}
Let $j\geq 2$ and $-j+\frac{1}{2}\leq s\leq-\frac{j}{2}$. Then
\begin{eqnarray}
      \left\|\Lambda^{-1}\partial_{x}(\prod_{j=1}^{2}u_{j})\right\|_{Z^{s}}\leq C\prod\limits_{j=1}^{2}\|u_{j}\|_{Z^{s}},
        \label{3.03}
\end{eqnarray}
\end{Lemma}
{\bf Proof.} Combining Lemmas 3.1, 3.2 with the definition of $Z^{s}$, we have that Lemma 3.3.

We have completed the proof of Lemma 3.3.

\bigskip
\bigskip

\noindent {\large\bf 4. Proof of Theorem  1.1}

\setcounter{equation}{0}

 \setcounter{Theorem}{0}

\setcounter{Lemma}{0}

\setcounter{section}{4}
Now we are in a position to prove Theorem 1.1.
We define
\begin{eqnarray*}
&&\Phi(u)=\eta(t) S(t)\phi-\frac{1}{2}\eta(t) \int_{0}^{t}S(t-t^{'})\eta (t^{'})\partial_{x}(u^{2})dt^{'},\label{4.01}\nonumber\\
&&B=\left\{u\in Z^{s}: \quad \|u\|_{ Z^{s}}\leq C\|\phi\|_{H^{s}(\mathbf{T} )}\right\}.
\end{eqnarray*}
By using   Lemmas 2.4, 2.6, 3.3,  we have that
\begin{eqnarray*}
&&\left\|\Phi(u)\right\|_{Z^{s}}\leq \left\| S(t)\phi\right\|_{Z^{s}}
+\left\|-\frac{1}{2}\eta(t) \int_{0}^{t}S(t-t^{'})\eta (t^{'})
\partial_{x}(u^{2})dt^{'}\right\|_{Z^{s}}\nonumber\\&&\leq C_{1}\|\phi\|_{H^{s}(\mathbf{T})}
+C\left\|\eta(t)\partial_{x}(u^{2})\right\|_{Z^{s}}\leq C\|\phi\|_{H^{s}(\mathbf{T})}+C\|u\|_{Z^{s}}^{2}.
\end{eqnarray*}
For $u,v\in B$,    provided that $\|\phi\|_{H^{s}(\mathbf{T})}$ is sufficiently small,  we derive that
\begin{eqnarray*}
&&\left\|\Phi(u)-\Phi(v)\right\|_{Z^{s}}\nonumber\\&&\leq C
\left(\|u\|_{Z^{s}}+\|v\|_{Z^{s}}\right)\|u-v\|_{Z^{s}}\nonumber\\
&&\leq 2C\|\phi\|_{H^{s}(\mathbf{T})}\|u-v\|_{Z^{s}}\leq \frac{1}{2}\|u-v\|_{Z^{s}}.
\end{eqnarray*}
For large initial data, if $u(x,t)$ is the solution to  (\ref{1.01})-(\ref{1.02}), then $u_{\mu}(x,t):=\mu^{-2j}u\left(\frac{x}{\mu},\frac{t}{\mu^{2j}}\right)$ is the solution to
\begin{eqnarray}
&&\partial_{t}u_{\mu}+(-1)^{j+1}\partial_{x}^{2j+1}u_{\mu}
     + \frac{1}{2}\partial_{x}\left[(u_{\mu})^{2}\right]
     = 0,\label{4.01}\\
    &&u_{\mu}(x,0)=\mu^{-2j}u_{0}\left(\frac{x}{\mu}\right),\quad x\in \mathbf{T}=\R/2\lambda\mu\pi, \label{4.02}
\end{eqnarray}
since $\|u_{\mu}(x,0)\|_{H^{s}}\leq C\mu ^{-2j+\frac{1}{2}-s}\|u_{0}\|_{H^{s}}$, we take $\mu$ sufficiently large, then
$\|u_{\lambda}(x,0)\|_{H^{s}}$  is sufficiently small, which is reduced to the case of  small initial data.

The proof of the rest of Theorem 1.1 can be found in  \cite{MT,TKato}, thus, we omit the process.

We have completed the proof of Theorem 1.1.

\bigskip
\bigskip

\noindent {\large\bf 5. Proof of Theorem  1.2}

\setcounter{equation}{0}

 \setcounter{Theorem}{0}

\setcounter{Lemma}{0}

\setcounter{section}{5}

This section is devoed to presenting Theorem 1.2.  Following the method of \cite{Bourgain97},  it suffices to derive that
\begin{eqnarray}
\left\|A_{3}(u_{0})\right\|_{\dot{H}^{s}}\leq C\|u_{0}\|_{\dot{H}^{s}}^{3}\label{6.01}
\end{eqnarray}
fails when $s<-j+\frac{1}{2},$ with $j\geq 2,j\in Z.$
where
\begin{eqnarray*}
&&A_{3}(u_{0})=2\int_{0}^{t}S(t-s) \partial_{x}(u_{1}(s)A_{2}(u_{0})(s),\nonumber\\
&&u_{1}(t)=S(t)u_{0},\nonumber\\
&&A_{2}(u_{0})=\int_{0}^{t}S(t-s)\partial_{x}\left[(u_{1}(s))^{2}\right]ds.
\end{eqnarray*}
Let
\begin{eqnarray*}
u_{0}=\phi_{N}=N^{-s}\left(\chi_{N}(k)+\chi_{-N}(k)\right).
\end{eqnarray*}
It is easily checked that
$\|u_{0}\|_{\dot{H}^{s}}\sim 1$. By a direct  computation,  we derive that
\begin{eqnarray}
\mathscr{F}_{x}A_{2}(u_{0})(t)=\sum_{k_{1}\neq 0,k\neq k_{1}}k\frac{e^{itp(k)}
-e^{itp(k_{1})+itp(k-k_{1})}}{q_{0}(k_{1},k-k_{1})}\mathscr{F}_{x}u_{0}(k_{1})
\mathscr{F}_{x}u_{0}(k-k_{1})\label{6.02},
\end{eqnarray}
where
\begin{eqnarray*}
q_{0}(k_{1},k-k_{1})=k_{1}^{2j+1}+k_{2}^{2j+1}-k^{2j+1},
\end{eqnarray*}
\begin{eqnarray*}
&&A_{3}(u_{0})=2\sum_{k_{1}\neq 0}\sum_{k_{2}\neq 0}\sum_{k_{3}\neq 0}e^{(k_{1}+k_{2}+k_{3})x+
ip(k_{1}+k_{2}+k_{3})t}\left(-\frac{1-e^{-i q_{1}t}}{q_{1}}+\frac{1-e^{-i q_{2}t}}{q_{2}}\right)\nonumber\\
&&\qquad \times \frac{(k_{1}+k_{2}+k_{3})(k_{2}+k_{3})}{q_{0}(k_{2},k_{3})}\prod_{j=1}^{3}\mathscr{F}_{x}u_{0}(k_{j}),
\end{eqnarray*}
and
\begin{eqnarray*}
&&q_{1}=k_{1}^{2j+1}+k_{2}^{2j+1}+k_{3}^{2j+1}-k^{2j+1},\nonumber\\
&&q_{2}=k_{1}^{2j+1}+(k_{2}+k_{3})^{2j+1}-k^{2j+1}.
\end{eqnarray*}
Obviously,
when $k_{1}=-N$ and $k_{2}=k_{3}=N,$ $q_{2}$ does not vanish but $q_{1}$ vanishes.
\begin{eqnarray}
 CN^{-2s-(2j-1)}\leq \left\|A_{3}(u_{0})\right\|_{\dot{H}^{s}}\sim 1.\label{6.03}
\end{eqnarray}
When $s<-j+\frac{1}{2}$, letting $N\rightarrow +\infty$ yields that the left hand side of (\ref{6.03}) goes to $+\infty$.
Thus, we obtain the contradiction.

We have completed the proof of Theorem 1.2.

\leftline{\large \bf Acknowledgments}

\bigskip

\noindent

 This work is supported by the Natural Science Foundation of China
 under grant numbers 11171116 and 11401180  and 11371367. The first author is also
 supported by   the Young core Teachers program of Henan Normal University and  15A110033.
  
\end{document}